\documentclass[12pt,a4paper]{article}

\usepackage{amssymb}
\usepackage{float}
\usepackage{amsthm}
\usepackage{lipsum}
\usepackage{comment}
\usepackage{amsfonts}
\usepackage{hyperref}
\usepackage{subcaption}
\usepackage{graphicx}
\usepackage{epstopdf}
\usepackage{algorithmic}
\usepackage{ulem}
\usepackage{xcolor}
\usepackage{mathtools}
\usepackage{etoolbox}

\newtheorem{theorem}{Theorem}
\newtheorem{definition}{Definition}
\newtheorem{lemma}{Lemma}
\newtheorem{proposition}{Proposition}
\newtheorem{corollary}{Corollary}
\newtheorem{remark}{Remark}
\newtheorem{example}{Example}

\usepackage[ruled,vlined]{algorithm2e}

\newcommand{\peru}{\text{\normalfont Per}_1}
\newcommand{\per}{\text{\normalfont Per}}

\newcommand{\pcr}{{\rm PCR}}
\newcommand{\bv}{{BV}}

\newcommand{\rofa}{{{\rm AROF}}}

\newcommand{\acv}{{\rm ACV}}
\newcommand{\cv}{{\rm CV}}
\DeclareMathOperator*{\argmin}{arg\,min}

\newcommand{\leb}{\mathcal{L}}

\newcommand{\z}{\mathbf{z}}

\newcommand{\Fcal}{\mathcal{F}}

\newcommand{\Div}{{\rm div}}

\graphicspath{ {./Figures/} }

\newcommand{\sm}[1]{{\color{black}{#1}}}
\newcommand{\vp}{\color{black}}

\newcommand{\Gcal}{\mathcal{G}}

\usepackage[hang,flushmargin]{footmisc}

\usepackage{titlesec}

\titleformat*{\section}{\normalsize\bfseries}
\titleformat*{\subsection}{\normalsize\bfseries}
\titleformat*{\subsubsection}{\normalsize\bfseries}
\titleformat*{\paragraph}{\normalsize\bfseries}
\titleformat*{\subparagraph}{\normalsize\bfseries}

\newcommand\rest[2]{{
		\left.\kern-\nulldelimiterspace 
		#1 
		\right|_{#2} 
}}

\begin{document}
	
\vspace{2ex}
\begin{center}
	{\Large\textbf{
			Anisotropic Chan--Vese segmentation
	}}
\end{center}
\begin{center}
{\large
		Salvador Moll$^{a, \ast}$, Vicent Pallardó--Julià$^{a,b}$\\
}
\end{center}

\vspace{8ex}
\begin{center}
\centering\begin{minipage}{\linewidth-1.5cm}
\noindent\rule{\linewidth}{1pt}

{\footnotesize
\noindent {\bf Abstract.} In this paper we study a variant to \textit{Chan--Vese} (CV) segmentation model with rectilinear anisotropy. We show existence of minimizers in the $2$-phases case and how they are related to the (anisotropic) \textit{Rudin-Osher-Fatemi} (ROF) denoising model. Our analysis shows that in the natural case of a piecewise constant on rectangles image ($\pcr$ function in short), there {exists} a minimizer of the CV functional which is also piecewise constant on rectangles over the same grid that the one defined by the original image. In the multiphase case, we show that minimizers of the CV multiphase functional also share this property in the case that the initial image is a $\pcr$ function. We also investigate a multiphase and anisotropic version of the Truncated ROF algorithm, and we {compare} the solutions given by this algorithm with minimizers of the multiphase anisotropic CV functional.

}
\end{minipage}
\end{center}
\let\thefootnote\relax\footnotetext{
	\noindent$^a$ { Department d'Anàlisi Matemàtica, Universitat de València, C/Dr. Moliner, 50, Burjassot, Spain}\\
	$^b$ {Kimera Technologies, Lanzadera, C/ del Moll de la Duana, s/n, València, Spain}\\
{	 $\empty^\ast$ Corresponding author: j.salvador.moll@uv.es}}

\section{Introduction}
Image segmentation consists in partitioning a given image into multiple segments in which pixels share some characteristics. One of the most relevant models in the field of image segmentation is the \textit{Mumford-Shah} (MS) model, introduced by the authors in \cite{MS}. This model was the seed of a very successful approach to the problem: variational techniques with level set formulations. A particular case of the Mumford-Shah model is the case in which the objective function is piecewise constant inside some domains with finite perimeters. This model was introduced by Chan and Vese in the $2$-phases and multiphases framework  in \cite{CV} and \cite{CV-n}, respectively. They are known as \textit{Chan-Vese} (CV) models and {play a cornerstone role in diverse recent applications of image processing (e.g. see \cite{CORONA, NAG, Zheng2022, ZHENG}).}\\

The starting point of this paper is the recent study \cite{Cai} about the relationship between the CV model in image segmentation and \textit{Rudin-Osher-Fatemi's} (ROF) model in image denoising (see \cite{ROF}). In \cite{Cai}, the authors show that a thresholding in the ROF model's solution provides a partial minimizer of the two phases CV functional, which can be {written} as
$$
{\rm CV}_2(\Lambda,c_1,c_2)={\text{Per}(\Lambda;\Omega)} + \mu\left(\int_{\Lambda}(c_1-f)^2\,dx+\int_{\Omega\setminus \Lambda} (c_2-f)^2\,dx\right).
$$

A partial minimizer in the sense that one can obtain a minimizer in each of its three variables; i.e. letting $\Lambda_u:=\{x: u(x)\geq \frac{c_1+c_2}{2}\}$, with $u$ being the minimizer of the ROF functional {(see \cite[Proposition 2.6]{CHAMBOLLE-CASELLES})}, one has

\begin{equation*}
	\begin{cases}
		\begin{aligned}
			\Lambda_u&\in \argmin_{\Lambda:  \per (\Lambda)<+\infty} {\rm CV}_2(\Lambda,c_1,c_2) \\
			\left(\textstyle \frac{1}{|\Lambda|}\int_\Lambda f\,dx, \, {\textstyle \frac{1}{|\Omega\setminus \Lambda|}\int_{\Omega\setminus \Lambda} f\,dx} \right)&\in  \argmin_{(c_1, c_2)\in [0,1]^2} {\rm CV}_2 (\Lambda,c_1,c_2)
		\end{aligned}
	\end{cases}.
\end{equation*}

{Despite this, whether} $\left(\Lambda_u,\frac{1}{|\Lambda_u|}\int_{\Lambda_u} f\,dx,\frac{1}{|\Omega\setminus \Lambda_u|}\int_{\Omega\setminus \Lambda_u} f\,dx\right)$ is a true minimizer of ${\rm CV}_2$ remains as an open problem. On the other hand, as stated in \cite{Cai}, it is of interest to understand if this relationship is still valid in some variants of both CV and ROF models.\\

Our work focuses in {particular} anisotropic variants of these models. {{Our motivation comes from }some {well known} features of the anisotropic $\ell_1$ version of the ROF model, such as sharp recovery of edges (see \cite{CHEN, ESEDOGLU-OSHER}), exact computability (e.g. see \cite{BOYKOV,GOLDSTEIN-OSHER,MR3709886}) and an observed reduction of the staircasing effect detected in the isotropic version (see \cite{CHEN, NIKOLOVA}). Because of that, we establish the anisotropic case as the $\ell_1$ one, i.e. we replace in the CV models the usual total variation by the total variation with respect to $|\cdot|_1$, defined by $|\boldsymbol{ v}|_1:=|v_1|+|v_2|\,$ for  $\boldsymbol{v}=(v_1,v_2)\in\mathbb{R}^2$.\\
}

{{Our} main objectives are the next ones: First, we show that}
there is a global minimizer of the two phases anisotropic CV functional
$$
\acv_\mu(\Lambda,c_1,c_2):={\per_1(\Lambda;\Omega)}+\mu\left(\int_{\Lambda}(c_1-f)^2\,dx+\int_{\Omega\setminus\Lambda}(c_2-f)^2\,dx\right),
$$
whose first component is an upper level set of a minimizer to the anisotropic ROF functional {in $L^2(\Omega)$, which is defined as follows:}
\begin{equation*}
	\rofa_\lambda(u) := |Du|_1(\Omega)+\frac{\lambda}{2}\int_{\Omega} (u-f)^2dx.
\end{equation*}

In order to show this result, we need to assume that $\Omega$ is a rectangle and that the data considered belong to a suitable space. Namely, we will assume that $f$ is piecewise constant on rectangles (denoted by $\pcr$ and defined in Def. \ref{def:pcr}). We point out that this restriction is harmless from the point of view of applications. Our strategy is as follows: First of all, we  generalize the $\acv$ functional, which is defined only for sets of finite perimeter, to an energy functional defined on $L^2(\Omega)\times [0,1]^2$ as follows
\begin{equation*}
	{\Gcal}_\mu(u,c_1,c_2) := |Du|_1(\Omega) +\mu  \int_\Omega (u(c_1 - f)^2  + (1 - u)(c_2 - f)^2)dx + \int_\Omega \mathbb{I}_{[0, 1]}(u)dx,
\end{equation*}
where $\mathbb{I}_{[0,1]}$ is the indicator function on $[0,1]$. In relation to the above, we note that $\mathcal G_\mu(\cdot, c_1, c_2)$ and $\rofa_\lambda(\cdot)$ are solely finite on {the space of bounded variation functions}, $\bv(\Omega)$. {On the other hand,} the indicator function restricts the range of $\Gcal_\mu(\cdot,c_1, c_2)$ to $[0,1]$. Additionally, we remark that if $E$ {is a set with finite perimeter} and $u=\chi_E$, ${\Gcal}_\mu(u,c_1,c_2)=\acv_\mu(E,c_1,c_2)$.\\

On these functionals, we prove existence of a global minimizer {of} ${\Gcal}_\mu$ through a direct variational method. Then, we show that a truncation of the solution to the $\rofa$ functional (i.e. an upper level set of the solution) yields the first component of a minimizer to $\Gcal_\mu$. To {show} that this is indeed the case, we rely on the description of explicit solutions of $\rofa$ {obtained} in \cite{MR3709886} for any PCR datum and on the corresponding Euler-Lagrange equations to both functionals. {In doing so}, we {find} a global minimizer both of $\mathcal G$ and of $\acv$. All the results concerning the 2-phases case are worked out in Section \ref{sec:2-phases}.

\begin{remark}\label{rmk:MP1}
	
	In the 1-dimensional case $(\Omega = [a,b]\subset\mathbb{R})$, it was shown in \cite{MP1} that there is a minimizer to the {\rm CV} problem with first component satisfying the following two properties:
	\begin{itemize}
		\item[(a)] The boundary of the set belongs to a sole level set of the datum $f$ in a multivalued sense.
		\item[(b)] If $f$ is piecewise constant, then the boundary of the set is contained in the jump set of $f$.
	\end{itemize}
	
	The above properties are shown to be false in the anisotropic case {as shown in Example \ref{ex:break}} as well as {it is in the two--dimensional case of} the standard {\rm CV} model (see \cite[Remark 4.2]{MP1}). However, as a by-product of our result, we obtain that, in the case of a {\rm PCR} datum, the first component of a minimizer to the $\acv_\mu$ functional is also a {\rm PCR} function; i.e the minimizing set is a rectilinear polygon. Moreover, its essential boundary belongs to a grid generated by the datum itself. This last property permits to design a trivial algorithm to compute a minimizer of the $\acv_\mu$ functional.
	
\end{remark}

{In second place, we deal with the $n$-phases anisotropic CV model}. In this case, we {decide} to slightly modify the original (anisotropic) functional, which reads as
\begin{equation}\label{A-CV-multiphase}
	\acv^n_\mu(\boldsymbol{\Omega},{\bf c}):= \sum_{i=1}^n \left({\per_1(\Omega_i;\Omega)} + \mu \int_{\Omega_i} (c_i - f)^2\,dx \right),
\end{equation}
{with $\mu > 0$ , $\boldsymbol{\Omega}:=\{\Omega_i\}_{i=1}^n\in\mathcal P_n (\Omega)$, {where, by $\mathcal P_n(\Omega)$ we denote the set of all non empty $n$-parts {(disjoint)} partition of $\Omega$} and $\boldsymbol{c} :=\{c_i\}_{i=1}^n\in[0,1]^n$}. {The proposed modification consists in} not counting more than once the length of the possible overlaps of the boundaries of the partition. To do this, we propose the following energy functional:
\begin{equation}\label{g-functional}
	\Gcal_\mu^n(\boldsymbol{\Sigma},\boldsymbol{c}) = \sum_{i=1}^{n} \left( \peru(\Sigma_{i-1}; \sm{\rm int(\Sigma_{i})})  + \mu\int_{\Sigma_i\setminus\Sigma_{i-1}} (c_{i}-f)^2dx  \right),
\end{equation}
{with $\mu > 0$ and  $\boldsymbol{c} := \{c_i\}_{i=1}^{n}\in[0,1]^n$ and $\boldsymbol{\Sigma} := \{\Sigma_i\}_{i=0}^n\in \mathcal{P}^*_n(\Omega)$, where $\mathcal{P}^*_n(\Omega) := \{\{\Lambda_i\}_{i=0}^{n}  :   \emptyset=\Lambda_0\subseteq\Lambda_i\subseteq\Lambda_{i+1}\subseteq\Lambda_n = \Omega\}$.}

\medskip

We observe that defining $\Omega_i:=\Sigma_i\setminus \Sigma_{i-1}$, the unique difference between both functionals is that in $\acv_\mu^n$ some edges will count more than twice \sm{(in the case that an edge belongs to the boundary of more than two different upper level sets)}  in the length term while in $\Gcal_\mu^n$ they are  counted only once.

\medskip
As in the $2$-phases case, the existence of a minimizer $(\boldsymbol{\Sigma^*},{\bf c^*})$ follows from the direct method in calculus of variations (see Proposition \ref{prop_2}).
Moreover, one can assume that $0\leq c^*_{i+1}\leq c^*_i\leq 1$ for any $i=1,\ldots, n-1$. Next, we observe that
\begin{equation*}
	\boldsymbol{\Sigma^*} \in \argmin_{\boldsymbol{\Lambda}\in \mathcal{P}^*_n(\Omega)} \Gcal_\mu^n(\boldsymbol{\Lambda}, {\boldsymbol c^*}) \quad \text{{with}} \quad c_i^*=\frac{1}{|\Sigma_{i}^*\setminus \Sigma_{i-1}^*|}\int_{\Sigma_{i}^*\setminus \Sigma_{i-1}^*}f(x)\,dx.
\end{equation*}

In the case that $f\in \pcr(\Omega)$, by the tools developed in \cite{MR3709886}, we can show that $\Sigma^*_i$ is a rectilinear polygon whose boundaries lie on the grid generated by $f$. These results concerning the multiphase case are the core of Section \ref{sec:multiphase}.\\

{Thirdly, we discuss {about a possible} relationship between the minimizers of a variant of $\acv_\mu^n$  and the truncated $\rofa_\lambda$ functional, from a general point of view. For this purpose, we define the following general $\acv_\mu^n$ variant as
	\begin{equation}\label{eq:cvphi}
		\cv^n_{\varphi,\boldsymbol \mu}(\boldsymbol\Omega, \ \boldsymbol{c}) := \sum_{i=1}^n\left( \per_\varphi(\cup_{j=1}^i\Omega_i; \ \Omega) + \mu_i\int_{\Omega_i}(c_i - f)^2\, dx  \right)
	\end{equation}
	where $\boldsymbol \mu := \{\mu_i\}_{i=1}^n\in [0,+\infty)^n$, $\boldsymbol{\Omega}:=\{\Omega_i\}_{i=1}^n\in\mathcal{P}_n(\Omega)$ and $\boldsymbol{c} := \{c_i\}_{i=1}^n\in[0,1]^n$. Similarly, the truncated ROF functional for $n$ phases is defined as
	\begin{equation}\label{eq:trofphi}
		{\rm TROF}^n_{\varphi, \lambda}(\boldsymbol\Sigma, \boldsymbol\tau) := \sum_{i=1}^{{n-1}}\left(\per_\varphi(\Sigma_i; \ \Omega) + \lambda\int_{\Sigma_i}(\tau_i - f)\, dx \right)
	\end{equation}
	where $\lambda >0$, $\boldsymbol{\Sigma}:=\{\Sigma_i\}_{i=0}^{{n}}\in\mathcal{P}^*_{{n}}(\Omega)$ and $\boldsymbol{\tau} := \{\tau_i\}_{i=1}^{{n-1}}\in[0,1]^{n-1}$. {Here,}  these generalizations depend on a positively 1-homogeneous convex function $|\cdot|_\varphi$. Note that $CV_{\varphi,\boldsymbol\mu}^n = \Gcal_\mu^n$ if $|\cdot|_\varphi$ is the 1-norm and $\mu_i = \mu$ for every $i$. In Section \ref{sec:trof} we prove that, {in the 2-phases case, it is possible to obtain a relationship between the minimizers of $\acv^2_{\mu}$ and ${\rm TROF}^2_{1,\lambda}$. However, we prove that a similar relationship in the multiphase case cannot be true, in general. }
	Finally, we remark the relationship between ${\rm TROF}_{\varphi,\lambda}^n$ and $\rofa_\lambda$ in the anisotropic case. \\
	
	The paper finishes with some applications of our results. In Section \ref{sec:applications}, we show the strength of them with some examples on 2-phases and multiphase image segmentation, by comparing them with similar isotropic processes.}

\section{Preliminaries}
Before starting, we {prescribe some notations and we } provide a basic knowledge {on bounded variation functions and anisotropies.}
\subsection{Notations}
Throughout the paper, $\Omega\subset \mathbb{R}^N$ will denote an open bounded set with boundary $\partial\Omega$ and $\nu^\Omega$ will denote the outer unit exterior normal at a point on the boundary, when defined. We denote by \sm{$L^p(\Omega)$ ($1\leq p<+\infty$), $L^\infty(\Omega)$} and $\mathcal M(\Omega)^M$ the set of Lebesgue $p-$integrable functions in $\Omega$, \sm{the set of essentially bounded measurable functions in $\Omega$} and the set of finite vector Radon measures on $\Omega$, respectively. In the case that $M=1$, we omit the index. {For any measurable set $E\subset\mathbb R^N$ with respect to the Lebesgue measure in $\mathbb R^N$, we denote by $|E|$ the Lebesgue measure of $E$.} By $\mathcal H^{N-1}$, we denote the {$(N - 1)-$dimensional Haussdorf} measure. We denote by $\bv(\Omega)$ the Banach space of functions of bounded variation in $\Omega$ with the norm defined next:
\begin{align*}
	\bv(\Omega)&:=\{u\in L^1(\Omega) :  Du\in\mathcal M(\Omega)^N\},\\
	\|u\|_{\bv(\Omega)}&:=\|u\|_{L^1(\Omega)}+\|D u\|_{\mathcal M(\Omega)^N}\,,
\end{align*}
where $Du$ is the distributional gradient of $u$. If $u\in \bv(\Omega)$, then the measure $Du$ can be decomposed into its absolutely continuous part and singular with respect to the Lebesgue measure $\mathcal L^N$: $$Du=\nabla u \mathcal L^N+D^s u\,,$$ with $\nabla u$ being the Radon-Nikodym derivative of $Du$ with respect to $\mathcal L^N$. We use standard notations and results for bounded variation functions such as those in \cite{ambrosio2000fbv}. In particular, $J_u$ will denote the jump set of $u\in \bv(\Omega)$ and $\nu^u$ the Radon-Nikodym of $Du$ with respect to its total variation; i.e. $\nu^u=\frac{Du}{|Du|}$. Moreover, $D^s u$ can be split into the jump part of the measure $D^j u$ and its Cantor part $D^cu$ and $$D^j u=(u^+-u^-)\nu^u\rest{\mathcal H^{N-1}}{J_u},$$ with $u^+$ and $u^-$ being the approximate limits at a jump point. We will assume the criterion that $\nu^u$ is oriented at $\mathcal H^{N-1}-$a.e. point on $J_u$ in such a way that $u^+>u^-$.

\medskip
We consider the space
\begin{equation*}
	X_\Omega=\left\{\z\in L^{\infty}(\Omega;\mathbb R^N)\,:\
	\,\Div \, \z\in L^2(\Omega)\right\}\,.
\end{equation*}
In \cite[Theorem 1.2]{anzellotti1}, the weak trace on the boundary of a
bounded Lipschitz domain $\Omega$ of the
normal component of $\z \in X_\Omega$ is defined. Namely, it is proved
that the formula
\begin{equation*}
	\left\langle [\z,\nu^\Omega]\,,\,\rho\right\rangle:=\int_{\Omega} \rho\,
	\Div \, \z + \int_{\Omega}\z\cdot\nabla\rho\,\qquad (\rho\in
	C^1\left(\overline{\Omega}\right)),
\end{equation*}
defines a linear operator $[\cdot,\nu^\Omega]: X_\Omega
\rightarrow
L^\infty(\partial \Omega)$ such that
\begin{equation*}
	\label{linftytrace}
	\Vert \, [\z,\nu^\Omega] \, \Vert_{L^\infty(\partial
		\Omega)} \leq \Vert \z
	\Vert_{L^\infty(\Omega)},
\end{equation*}
for all $\z \in X_\Omega$ and $[\z,\nu^\Omega]$ coincides with the
pointwise trace of the normal component if $\z$ is smooth. {From now on, the above definitions will be considered in the special case that $N = 2$.}

\subsection{Anisotropies and Chan-Vese model}
Next, we refer about \textit{isotropic} and \textit{anisotropic} models. This classification depends on the norm in the \textit{total variation} used in each model. As follows, we introduce its definition according to \cite{AmarBellettini}.
\begin{definition}\label{def:varphiTV}
	Let $\Omega\subseteq \mathbb{R}^2$ be an open set, $|\cdot|_\varphi:\mathbb R^2\to [0,+\infty[$ be a convex, positively 1-homogeneous function such that $|{\bf x}|_\varphi>0$ if ${\bf x}\neq 0$ and {let} $u$ be a function in $L^2(\Omega)$. Then, the {total variation of ${u}$ with respect to $\boldsymbol{|\cdot|_\varphi}$}, denoted by $|Du|_{\varphi}(\Omega)$, is defined as
	\begin{equation*}
		|Du|_{\varphi}(\Omega) := \sup \left\{  \displaystyle \int_\Omega u \,\Div \, \boldsymbol{\eta} \, dx \, : \, \boldsymbol{\eta}\in C_{c}^1 (\Omega; \mathbb{R}^{2}) , \, |\boldsymbol{\eta}|_{\varphi}^* \leq 1  \right\},
	\end{equation*}
	{where $|\cdot|_\varphi^*$ denotes the convex conjugate of $|\cdot|_\varphi$. We note that $|Du|_\varphi(\Omega) < +\infty$ if and only if $u\in \bv(\Omega)$, in the case that {$|\cdot|_\varphi$} is a norm. {In this context, we write $B_\varphi(x;R)$ the ball with radius $R>0$ centered at $x\in\mathbb R^N$ with respect to the $|\cdot|_\varphi$ distance.} In addition, 	we denote by $\per_\varphi(E;\Omega) = |D\chi_E|_\varphi(\Omega)$ the {${\boldsymbol{\varphi}}$-perimeter} of E in $\Omega$, where $E\subseteq\Omega$.}{ When $\Omega=\mathbb R^N$, we omit it and we simply write $\per_\varphi(E)$.} {Finally, for a finite perimeter set $E$, we denote by $\partial^*E$ its reduced boundary, according to the notation in \cite{bellettininovagapaolini2}.}
\end{definition}

\begin{remark}
	{If $|\cdot|_\varphi$ is the Euclidean norm, {the definition above coincides with the total variation of the measure $Du$.} The models {in which the Euclidean norm} is used are called \textit{isotropic}. Otherwise, we will speak about \textit{anisotropies}, the models where these are used are called \textit{anisotropic}. Throughout this work, if the function $|\cdot|_\varphi$ is not explicitly stated, it will refer to the anisotropic model when $|\cdot|_\varphi$ is the 1-norm, i.e. $|(v_1 , v_2)|_1:=|v_1|+|v_2|\,$ for  $(v_1,v_2)\in\mathbb{R}^2$; and the $\ell_1$ total variation of $u\in \bv(\Omega)$ will be denoted by $|Du|_{1}(\Omega)$.}
\end{remark}

Next we introduce the Chan-Vese model in a general $|\cdot|_\varphi$ setting:
\begin{definition}\label{def:ACVmodel}
	Let $\Omega\subset \mathbb{R}^2$ be an open set, $|\cdot|_\varphi$ be as Def. \ref{def:varphiTV} and $f$ be an $L^2(\Omega)$ function such that $f(\Omega)\subseteq [0,1]$. We consider the following functional:
	\begin{equation*}\label{eq:APCMSfunctional}
		\cv_{\varphi, \mu}(\Sigma, c_1, c_2) = \per_\varphi(\Sigma; \Omega) + \mu\left(\int_{\Sigma}(c_1 - f)^2dx + \int_{\Omega\setminus\Sigma}(c_2 - f)^2dx\right),
	\end{equation*}
	where $\Sigma\subseteq\Omega$, $c_1, c_2\in [0,1]$ and $\mu >0$. The associated model to the above functional consists in finding a 3-tuple $(\Sigma, c_1, c_2)$ such that
	\begin{equation*}
		(\Sigma, c_1, c_2)\in\argmin_{\substack{\Lambda\subseteq\Omega \\ s_1,s_2\in[0,1]}}\cv_{\varphi, \mu}(\Lambda, s_1, s_2) \qquad \text{for a fixed $\mu$.}
	\end{equation*}
	According to this definition, the $\cv$ and $\acv$ models correspond to $\cv_{\varphi, \mu}$ model when $|\cdot|_\varphi=|\cdot|_2$ and $|\cdot|_\varphi=|\cdot|_1$, respectively; and their functionals are denoted by $\cv_\mu$ and $\acv_\mu$, respectively.
\end{definition}

\subsection{PCR functions and ROF problems}
\begin{definition}
	Let $\Omega\subset \mathbb{R}^2$ be an open set, $|\cdot|_\varphi$ be as Def. \ref{def:varphiTV} and $f$ be an $L^2(\Omega)$ function. The $\varphi$-{\rm ROF}$_\lambda$ functional is defined as follows:
	\begin{equation}\label{eq:A-ROFproblem}
		{\rm ROF}_{\varphi,\lambda}(u) = {|D u|_\varphi(\Omega)} + \frac{\lambda}{2}\int_\Omega (u-f)^2dx,
	\end{equation}
	where $u\in L^2(\Omega)$ and $\lambda>0$. The associated model to the above functional consists in finding the minimum of ${\rm ROF}_{\varphi, \lambda}$, i.e., solving the problem
	\begin{equation*}
		\min_{u\in L^2(\Omega)} \,\, {\rm ROF}_{\varphi, \lambda}(u)
	\end{equation*}
	
	{The minimizer, which we denote by $w_\lambda$, satisfies the associated Euler-Lagrange equation to ${\rm ROF}_{\varphi, \lambda}$, which is the next one:
		\begin{equation}\label{eq:ELphiROF}
			-\Div (\boldsymbol{z}_{w_\lambda}) = \lambda(w_\lambda - f),
		\end{equation}
		with $\boldsymbol{z}_{w_\lambda} \in \partial |Dw_\lambda|_\varphi$ {(here the symbol $\partial$ denotes the subdifferential of the anisotropic total variation; {see \cite{Moll} for a characterization})}.   According to these definitions, we denote by ${\rm ROF}_\lambda$ and $\rofa_\lambda$ models those ${\rm ROF}_{\varphi, \lambda}$ ones where $|\cdot|_\varphi$ are the Euclidean norm and 1-norm, respectively.}
\end{definition}

As we work with the $\rofa$ model when $f$ belongs to a particular family of functions, $\pcr(\Omega)$, we next introduce the definition of this space.

\begin{definition}\label{def:pcr}
	Let $\Omega$ be a rectangle. We say that $w$ is {piecewise constant on rectangles} and we write $w\in \pcr(\Omega)$ if $w$ has a finite number of level sets of positive $\leb^2$ measure, and each one is a rectilinear polygon up to an $\leb^2$-null set.
\end{definition}

In order to work with the results introduced in \cite{MR3709886}, we need to present some notation used in that paper. These terms {are introduced now.}

\begin{definition}
	{Let $\Omega$ be a bounded {rectangle in} $\mathbb{R}^2$. We denote by $G$ {any} set of horizontal and vertical lines on $\Omega$, which we will call a {grid} on $\Omega$, and any of the rectangles in the partition generated by the grid we will call it {cell}. Let $F$ be a rectangular polygon contained in $\Omega$ and let $f\in \pcr(\Omega)$. We denote
		\begin{itemize}
			\item by $\mathcal{F}(G)$ the set of all rectangular polygons which satisfy that each one of their sides is a segment of adjacent vertices of $G$.
			\item by $G(F)$ the minimal grid such that each side of $F$ is contained in one line of $G$, by $\mathcal{Q}_f$ the partition of $\Omega$ provided by the level-sets of $f$ and by $G_f$ the grid $\cup_{\Sigma\in \mathcal{Q}_f} G(\Sigma)$.
		\end{itemize}
	}
\end{definition}

It is immediate to observe that $\pcr(\Omega)\subset \bv(\Omega)$. Moreover, $J_u=\bigcup_{F\in \mathcal Q_u} \partial F$ for $u\in \pcr(\Omega)$. In \cite{MR3709886}, the subdifferential of the anisotropic total variation (in the case that $|\cdot|_\varphi=|\cdot|_1$) on PCR functions was characterized as follows:

\begin{lemma}\label{lem: subdifferential}
	Let $f\in \pcr(\Omega)$. Then, $v\in \partial |Df|_1$ if and only if $v\in L^2(\Omega)$ and there exists $\z\in X_\Omega$ such that $v=-\Div\, \z$, $\|\z\|_\infty\leq 1$ a.e. and $$[\z,\nu^\Omega]=0\,,\ [\z,\nu^F](x)=\left\{\begin{array}
		{cc} -1 & {\rm if \ } \rest{f}{F}<\rest{f}{F'}\,  x\in\partial F\cap\partial F'  \\ 1 & {\rm if \ } \rest{f}{F}>\rest{f}{F'} \, x\in\partial F\cap\partial F'
	\end{array}\right.,\mathcal H^{1}-{\rm a.e. \ on \ }\partial F\,,$$ with $F\neq F'\in \mathcal Q_f$.
\end{lemma}

Given $f\in \pcr(\Omega)$, in \cite[Theorem 5]{MR3709886} a minimizer $w_\lambda$ of {$\rofa_\lambda$} was obtained via a finite algorithm. In particular, it can be constructed as a $\pcr(\Omega)$ function with all level sets $\{F_i\}_{i=1}^l\subseteq \mathcal F(G_f)$. Moreover, $w_\lambda|_{F_k}>w_\lambda|_{F_{k+1}}$,
$$
\peru(\{w_\lambda>\tau\};\sm{\Omega})+{\lambda}\int_{\{w_\lambda>\tau\}}(w_\lambda-f)\,dx=0\,,\quad  {\rm for \ any \ }0<\tau<1,
$$ and, by $\eqref{eq:ELphiROF}$ {and Lemma \ref{lem: subdifferential}}, we obtain that there exists $\z_{w_\lambda}\in\partial |Dw_\lambda|_1$ such that
$$
\left\{\begin{array}{cc} -{\rm div} (\z_{w_\lambda})=\lambda(w_\lambda-f)  & \\
	\displaystyle [\z_{w_\lambda},\nu^{\{w_\lambda>\tau\}}]=-1 &  {\rm for \ any \ } 0<\tau<1.
\end{array}\right.
$$

Observe that, thanks to Lemma \ref{lem: subdifferential}, we obtain that
\begin{equation*}
	\label{subdifferential_inclusion}\partial |D w_\lambda|_1\subseteq \partial |D{{\chi_{\{w_\lambda>\tau\}}}}|_1\,, \quad {\rm for \ any \ }0<\tau<1.
\end{equation*}

\section{2-phases ACV model}\label{sec:2-phases}

{In this section we prove the existence of a $\pcr$ minimizer of the $2$-phases anisotropic Chan-Vese problem and its relationship with solutions to the $\rofa$ problem.}

\begin{definition}
	Let $\Omega\subset\mathbb{R}^2$ {be an open set} and let $f\in L^2(\Omega)$. Then, we define the generalized functional of $\acv_\mu$ as follows: ${\Gcal_\mu}:L^2(\Omega)\times [0,1]^2\to [0,+\infty[$,
	\begin{equation*}
		\begin{aligned}\label{eq:general-problem}
			{\Gcal_\mu}(u, c_1, c_2)=  & \, |Du|_1(\Omega) + \mu \int_\Omega \left(u(c_1 - f)^2 + (1 - u)(c_2 - f)^2\right) dx + \int_\Omega \mathbb{I}_{[0, 1]}(u),
		\end{aligned}
	\end{equation*}
	where $\mu > 0$. Concerning the generalised functional proposed above, we note that the partition $\{\Lambda, \Omega\setminus\Lambda\}$ has been relaxed by $u$ and $1-u$, providing a convex behaviour in terms of $u$. This function {$u$} has a range restriction on $[0,1]$, imposed by the indicator function of the interval $[0,1]$;
	\begin{align*}
		\mathbb{I}_{[0,1]}(x) &:= \begin{cases}
			0 \quad &\text{if $x\in[0,1]$}
			\\  +\infty &\text{otherwise.}
		\end{cases}
	\end{align*}
	
	The associated model to the above functional consists in finding a 3-tuple $(u, c_1, c_2)$ such that
	\begin{equation*}
		(u, c_1, c_2) \in \argmin_{\substack{ {w} \in L^2(\Omega)\\ (a,b)\in[0,1]^2}} \Gcal_\mu(w,a,b).
	\end{equation*}
\end{definition}

\begin{proposition}\label{prop:existence}
	There exists $(u,c_1,c_2)\in \bv(\Omega)\times[0,1]^2$ such that it is miminizer of $\Gcal_\mu$ and  $u(\Omega)\subseteq [0,1]$.
\end{proposition}
\begin{proof}
	The proof easily follows from the direct method in the calculus of variations. We give a sketch of it by the sake of completeness. Since $\Gcal_\mu\geq 0$ in $L^2(\Omega)\times [0,1]$, there exists $m=\inf \{\Gcal_\mu({\bf w}): {\bf w}\in L^2(\Omega)\times [0,1]^{\sm{2}}\}$. We take a minimizing sequence $\{(u_n,c_{1,n},c_{2,n})\}_{n\in \mathbb N}$. Since $u_n(\Omega) \subseteq [0,1]$, then there is a subsequence (not relabeled) and $u\in L^2(\Omega)$ such that $u_n\rightharpoonup u$ weakly in $L^2(\Omega)$ and that $u(x)\in [0,1]$ a.e. in $\Omega$. Moreover, from the lower semicontinuity of the {anisotropic} total variation and the coercivity of the $|\cdot|_1$ norm, we may assume that $u\in \bv(\Omega)$ and that $u_n\to u$ in $L^1(\Omega)$. On the other hand, since $c_{i,n}\in [0,1]$, for $i=1,2$, then there is a subsequence (not relabeled) and $c_i\in [0,1]$ such that $c_{i,n}\to c_i$ for $i=1,2$. Finally, from the lower semicontinuity of the anisotropic total variation we get that $$
	\Gcal_\mu(u,c_1,c_2)\leq \liminf  \Gcal_\mu(u_n,c_{1,n},c_{2,n})=m\,, $$
	which shows that $(u,c_1,c_2)$ is a minimizer.
\end{proof}

We now fix $0\leq c_2<c_1\leq 1$ and consider $\Gcal_{\mu,c_1,c_2}: L^2(\Omega)\to [0,+\infty]$ defined by
\begin{equation}\label{restricted functional}
	\Gcal_{\mu,c_1,c_2}(u):=\Gcal_\mu(u,c_1,c_2).
\end{equation}
{Observe that, in the case that $c_1=c_2$, the only minimizers are \sm{constant functions.}}

\begin{theorem}\label{thm:ROF-F}
	Let $0\leq c_2<c_1\leq 1$ and let $w$ be the minimizer of \eqref{eq:A-ROFproblem} with parameter $\lambda=\mu(c_1-c_2)$. Then,
	\[
	{u} =  \chi_{\Sigma} \quad \text{s.t.} \quad \Sigma:= {\left\{ x\in\Omega : w(x){>}\frac{c_1 + c_2}{2} \right\}}
	\]
	is a minimizer of $\Gcal_{\mu,c_1,c_2}$.
\end{theorem}
\begin{proof}
	First of all, we observe that the functional $\Gcal_{\mu,c_1,c_2}$ is convex and lower semicontinuous in $L^2(\Omega)$. Therefore, the existence of a minimizer $u\in \bv(\Omega)$ to \eqref{restricted functional} is guaranteed. Moreover, {applying \cite[Theorem 2.1]{MP1}}, its Euler-Lagrange equation is given by
	\begin{equation}\label{eq:EL_CVA}
		\Div \boldsymbol{z_{u}} = \mu(c_1 - c_2)(c_1 + c_2 - 2f) + h\,,
	\end{equation}
	where $\boldsymbol{z_{u}}\in \partial|D u|_1$ and $h\in \partial \mathbb{I}_{[0, 1]}(u)$.
	
	On the other hand,
	as $w$ is minimizer of \eqref{eq:A-ROFproblem}, the Euler-Lagrange equation with respect to this functional is fulfilled, i.e., we have
	\begin{equation}\label{eq:EL_ROF}
		{\Div\, z_{w}} = {\mu(c_1-c_2)(w - f)}\,,
	\end{equation}
	where $\boldsymbol{z_{w}}\in \partial|D w|_1\subset \partial |D u|_1$ by \eqref{subdifferential_inclusion}. Thus, we let $\boldsymbol{z_u} := \boldsymbol{z_{w}}\in \partial |Du|_1$  and we will prove that \eqref{eq:EL_CVA} is satisfied by $\boldsymbol{z_w}$. Proving this, we show that $u$ is a {minimizer of $\Gcal_{\mu,c_1,c_2}$}.
	
	From \eqref{eq:EL_ROF}, we get that \eqref{eq:EL_CVA} is satisfied if, and only if
	\begin{equation*}
		{\mu(c_1-c_2) \left(w - \frac{c_1 + c_2}{2}\right)} = h.
	\end{equation*}
	
	Therefore, we only need to show that $h\in \partial \mathbb{I}_{[0, 1]}(u)$. We distinguish two cases:
	\begin{itemize}
		\item If $w < \displaystyle\frac{c_1 + c_2}{2}$, we have that $u = 0$ and $h\in ]-\infty, 0] = \partial {\mathbb I}_{[0, 1]}(0)$.
		\item If $w \geq \displaystyle\frac{c_1 + c_2}{2}$, we have that $u = 1$ and $h\in [0, +\infty[ \, = \partial {\mathbb I}_{[0, 1]}(1)$.
	\end{itemize}
	Therefore, $u$ satisfies \eqref{eq:EL_CVA} and the proof is finished.
\end{proof}

\begin{remark}\label{rm:c1c2constant}
	Let $w$ be a minimizer of \eqref{eq:A-ROFproblem}. \sm{Since 
		each coordinate of a minimizer of $\Gcal_\mu$ fulfills its Euler-Lagrange equation, if}
	\begin{equation}\label{eq:c1c2nodepend}
		c_1 = \frac{\int_\Sigma f \, dx}{|\Sigma|}, \quad c_2 = \frac{\int_{\Omega\setminus \Sigma}}{|\Omega \setminus \Sigma|} \quad \text{where $\Sigma := \left\{ x \in \Omega \, : \, w(x)\geq \frac{c_1 + c_2}{2} \right\}$}\,,
	\end{equation}\sm{then, $(\chi_\Sigma,c_1,c_2)$ is a minimizer of $\Gcal_\mu$ if $\mu=\lambda/(c_1-c_2)$.}
	The existence of such a triplet $(\Sigma, c_1, c_2)$ satisfying \eqref{eq:c1c2nodepend} is guaranteed by Proposition 1 and Theorem 1.
\end{remark}

Now we conclude with \sm{the existence and }characterization of \sm{some} minimizers of $\acv_\mu$:

\begin{theorem}\label{thm:CV1sol}
	Let $w$ be the minimizer of $\rofa_\lambda$. \sm{Then, there is $0\leq c_2^*<c_1^*\leq 1$ 
			\sm{verifying} \eqref{eq:c1c2nodepend}, $(\chi_\Sigma, c_1^*, c_2^*)$ is a minimizer of $\mathcal G_\mu$ and 
				$(\Sigma, c_1^*, c_2^*)$ is a minimizer of $\acv_\mu$ if $\mu = {\lambda}/{(c_1^* - c_2^*)}$.}
		\end{theorem}
		\begin{proof}
			First, {we note} that it is obvious that, given $E\subset \Omega$ of finite perimeter,
			\[
			\mathcal G_\mu(\chi_E, c_1, c_2) = \acv_\mu(E, c_1, c_2).
			\]
			Then,
			\begin{equation}\label{eq:ineqF-ACV}
				\min_{\substack{u\in L^2(\Omega)\\ c_1, c_2\in[0,1]}} \mathcal G_\mu(u, c_1, c_2) \leq \min_{\substack{\Lambda\subset\Omega\\ c_1, c_2\in[0,1]}} \acv_\mu (\Lambda, c_1, c_2)
			\end{equation}
			On the other hand, we know that there exist minimizers of $\Gcal_\mu$ by Proposition \ref{prop:existence} and each coordinate of each minimizer satisfies its respective Euler-Lagrange equation. Let us denote by $(u^*,c_1^*,c_2^*)$ one minimizer. By Theorem \ref{thm:ROF-F}, $\Sigma$ satisfies \begin{equation*}
				\label{estimate_one_side}\sm{\Gcal}_\mu(u^*,c_1^*,c_2^*)=\sm{\Gcal}_\mu({\vp\chi_\Sigma},c_1^*,c_2^*)\geq \min_{ c_1, c_2\in[0,1]} \acv_\mu (\chi_{\Sigma}, c_1, c_2).
			\end{equation*}
			
			Therefore, $(c_1^*,c_2^*)$ fulfill the equalities in Remark \ref{rm:c1c2constant}. Finally, \eqref{eq:ineqF-ACV} and \eqref{estimate_one_side} imply that $$
			\min_{\substack{\Lambda\subset\Omega\\ c_1, c_2\in[0,1]}} \acv_\mu (\Lambda, c_1, c_2)\leq \sm{\Gcal}_\mu({\vp\chi_\Sigma},c_1^*,c_2^*) \leq \min_{\substack{\Lambda\subset\Omega\\ c_1, c_2\in[0,1]}} \acv_\mu (\Lambda, c_1, c_2),
			$$
			thus showing that $(\Sigma,c_1^*,c_2^*)$ is a minimizer of $ACV_\mu$.
		\end{proof}
		
		{\begin{example}
				\label{ex:break}
				Next, we show that the properties described in Remark \ref{rmk:MP1}, are not satisfied in the anisotropic 2-dimensional case either, using the previous result and the characterization of $\rofa_\lambda$'s solutions in \cite{MR3709886}. For that, we provide a counterexample.
				
				Let us consider $\Omega = [-1,1]^2$ and $ f= \chi_{A_1\cup A_2}$ \sm{with}
				\begin{equation*}
					A_1 = [-1/2,\ 1/2]^2, \quad A_2=[1/4, \ 1/2]\times[1/2,\ 3/4].
				\end{equation*}
				According to \cite[Theorem 5]{MR3709886}, the $AROF_{{16}}$ minimizer in this case is
				\begin{equation*}
					w = \frac{3}{4}\chi_{A_1} + \frac{1}{2}\chi_{A_2} + \frac{9}{{94}}\chi_{\Omega\setminus(A_1\cup A_2)}\,.
				\end{equation*}
				
				By Theorem \ref{thm:CV1sol}, we know that $(A_1,\ 1, \ 1/48)$ is a minimizer of $\acv_{{768}/47}$, whose jump set is not contained in the jump set of $f$ (and, therefore, in a sole level set of $f$).
				
		\end{example}}
		
		Finally, in algorithm \ref{alg:ACVsketch} we present the classical alternate algorithm that will allow us to approximate a solution of the $\acv$ model, depending on a tolerance term $\varepsilon_{tol}$ and an iteration maximum $n_{max}$.
		
		\medskip
		
		\begin{algorithm}[H]
			\caption{$\acv$ \sm{approximate} minimizer}
			\label{alg:ACVsketch}
			\begin{algorithmic}
				\STATE{\textbf{Initiation:} $k = 0$, $\mu>0$, $(c_1, c_2)\in[0,1]$ s.t. $c_1 > c_2$.}
				\STATE{$w_0 \gets$ $\rofa_{\mu(c_1 - c_2)}$ minimizer.}
				\STATE{$\Lambda \gets$ $\acv_\mu$ minimizer via Theorem {\ref{thm:ROF-F}}.}
				\vspace{0.2cm}
				\WHILE{$|\Lambda_k \bigtriangleup \Lambda_{k-1}|^2 > \varepsilon_{tol}\, \land \, k < n_{max}$}
				\vspace{0.2cm}
				\STATE{$c_1 \gets \frac{1}{|\Lambda_k|}\int_{\Lambda_k} f\,dx$, \, $c_2 \gets \frac{1}{|\Omega\setminus\Lambda_k|} \int_{\Omega\setminus\Lambda_k} f\,dx$, \quad $k \gets k + 1$}
				\vspace{0.1cm}
				\STATE{$w_k \gets$ $\rofa_{\mu(c_1 - c_2)}$ minimizer}
				\STATE{$\Lambda_k \gets$ $\acv_\mu$ minimizer via Theorem {\ref{thm:ROF-F}}}
				\vspace{0.2cm}
				\ENDWHILE
				\RETURN{$(\Lambda_k,c_1,c_2)$}
			\end{algorithmic}
		\end{algorithm}
		
		\medskip
		
		{Observe that we are imposing the condition $c_1>c_2$ at each step. This condition is guaranteed if initially it is so. Since the proof of this fact is a special case of the multiphase case, we refer to Proposition \ref{prop3}.}

		\section{Multiphase ACV model}\label{sec:multiphase}
		
		{This section is devoted to the study of the anisotropic multiphase Chan-Vese model. First of all, we prove that there exists a minimizer to \eqref{g-functional}.
			
			\begin{proposition}\label{prop_2}
				Let $\Omega\subset\mathbb{R}^2$ be an open set and $f\in L^2(\Omega)$. Then, there exists $\boldsymbol{\Sigma}\in \mathcal P^*_n(\Omega)$ and ${\bf c}\in [0,1]^n$ such that $$(\boldsymbol{\Sigma},{\bf c})\in \argmin_{\substack{\boldsymbol{\Lambda}\in \mathcal{P}^*_n(\Omega)\\ {\bf a}\in [0,1]^n}} \Gcal_\mu^n(\boldsymbol{\Lambda}, {\boldsymbol a}). $$
			\end{proposition}
			\begin{proof}
				The proof is very similar to the proof of Proposition \ref{prop:existence}. In particular, for a minimizing sequence $(\boldsymbol{\Sigma}^k,{\bf c^k})={((\Sigma^k_0,\ldots,\Sigma_n^k);(c^k_1,\ldots,c_n^k))}$, we can take $u^k_i:=\chi_{\Sigma_i^k}$ and work exactly as in Proposition \ref{prop:existence}. Finally, from the convergence a.e. obtained from the lower semicontinuity of the anisotropic total variation, we conclude that the weak limits in $L^2(\Omega)$, $u_i$ are of the form $u_i=\chi_{\Sigma_i}$ and that $\boldsymbol{\Sigma}:=(\Sigma_0,\ldots,\Sigma_n)\in\mathcal P^*_n(\Omega)$. The rest of the proof is totally analogous and we omit it.
		\end{proof}}
		
		\begin{theorem}\label{PCR}
			If ${\bf c}=\{c_i\}_{i=1}^n$ satisfies $1\geq c_i\geq c_{i+1}\geq 0$, then there exists a minimizer $\boldsymbol{\Sigma}$ to $\mathcal{G}^n_\mu(\cdot, \boldsymbol{c})$ for any initial datum $f\in \pcr(\Omega)$ and it satisfies that each one of its components {belongs to} $\Fcal(G_f)$; i.e: $\boldsymbol{\Sigma}$ is a rectangular partition of $\Omega$ whose components have their boundary on $G_f$.
		\end{theorem}
		
		\begin{proof}
			First of all, we rewrite {$\Gcal^n_\mu(\cdot,{\bf c})$} as follows:
			\begin{equation*}
				\Gcal^n(\boldsymbol{\Sigma}):={\Gcal^n_\mu(\boldsymbol{\Sigma},{\bf c})}  = \sum_{i=1}^{n-1} \left( \peru(\Sigma_i;\sm{{\rm int} (\Sigma_{i+1})}) + \mu\int_{\Sigma_i} g_{c_i, c_{i+1}}\,dx \right).
			\end{equation*}
			where $g_{c_{i}, c_{i+1}}(x) = (c_i - c_{i+1})(c_i + c_{i+1} - 2f(x))$.
			
			\bigskip
			{We now proceed as in \cite[Lemma 2]{MR3709886}. Given $\boldsymbol{\Sigma}\in \mathcal P_n^*(\Omega)$, and $\varepsilon>0$, we will construct $\boldsymbol{\Sigma^*}\in \mathcal P_n^*(\Omega)$ such that $\Sigma_i^*\in \mathcal F(G_f)$ and $$
				\Gcal^n(\boldsymbol{\Sigma^*})<\Gcal^n(\boldsymbol{\Sigma})+\varepsilon.$$
				We divide the proof into three steps:}
			\begin{enumerate}
				\item \textit{Smoothing}: We will construct a variant of $\boldsymbol{\Sigma}$, denoted by $\boldsymbol{\tilde{\Sigma}}$ such that $\boldsymbol{\tilde{\Sigma}}\in \mathcal P_n^*(\Omega)$ and their components have smooth boundaries. For that, we define $\tilde{\Sigma}_{i, \delta}$ as
				$$\tilde{\Sigma}_{i, \delta} = \{x\in\Omega : \psi_{\delta}\ast \chi_{\Sigma_{i}} {\geq} t\}\,,$$
				where $\psi_{\delta}$ is a standard smooth approximation of unity, $\delta$ is a positive parameter close to $0$ and $t\in (0,\frac{1}{2})$. As a consequence of this definition and following the idea of {\cite[Lemma 2]{MR3709886}, we know that $\tilde{\Sigma}_{i, \delta}$ is smooth and
					it is possible to set values $t$ and $\delta$ such that
					$$\Gcal^n(\boldsymbol{\tilde\Sigma})< \Gcal^n(\boldsymbol{\Sigma})+\varepsilon,$$}
				with $\boldsymbol{\tilde\Sigma}=\{\tilde\Sigma_i\}_{i=1}^{n-1}{\in \mathcal P_n^*(\Omega)}$.
				
				{\item \textit{Squaring}: Next, we will construct a variant of $\boldsymbol{\tilde{\Sigma}}$, denoted by $\boldsymbol{\hat{\Sigma}}$, such that {$\boldsymbol{\hat{\Sigma}}\in \mathcal P_n^*(\Omega)$} and their components have  rectangular boundaries. For this end, we will apply a squaring process on the components of $\boldsymbol{\tilde{\Sigma}}$ in each cell of $G_f$, avoiding {an} increase of {the} energy at $\Gcal^n$ . So, we consider a closed cell $C$ of $G_f$ and we suppose that $C^\circ\cap (\cup_{i=1}^{n-1}\partial \tilde \Sigma_i)=C^\circ\cap (\cup_{i=m_1}^{m_2}\partial \tilde \Sigma_i)$. {We split $\{\Sigma_i\}_{i=m_1}^{m_2}$} into those sets whose index satisfies $g_{c_i, c_{i+1}}>0$ in the cell $C$ and the rest, denoting them by $\{\tilde{\Sigma}_i\}_{i=m_1}^{i_0-1}$ and $\{\tilde{\Sigma}_i\}_{i=i_0}^{m_2}$ respectively, providing a well defined division thanks to the condition $c_i\geq c_{i+1}$, {and the fact that $f$ is constant in $C^\circ$}.
					
					First, we work with the sets $\tilde{\Sigma}_i\in\{\tilde{\Sigma}_i\}_{i=i_0}^{m_2}$. Following the idea of \cite[Lemma 2, Step 2]{MR3709886}, we construct a covering $W^{i}$ of $\partial (\tilde{\Sigma}_i\cap C^\circ)$, made up of squares, such that
					\begin{equation}
						\peru(\tilde{\Sigma}_i\cup{W^{i}};\sm{{\rm int} (\tilde \Sigma_{i+1})}) \leq \peru(\tilde{\Sigma}_i;\sm{{\rm int} (\tilde \Sigma_{i+1})}) \, \text{and} \, \leb^2(\tilde{\Sigma}_i\cup{W^{i}}) \geq \leb^2(\tilde{\Sigma}_i)
						\label{eq:nonincrease}
					\end{equation}

					Now, we apply a similar argument on $\{\tilde\Sigma_i\}_{i=m_1}^{i_0 - 1}$. In this case, we repeat the previous approach on $\Sigma_{i+1}\setminus\tilde\Sigma_i\in\{\tilde\Sigma_{i+1}\setminus\tilde\Sigma_i\}_{i=m_1}^{i_0 - 1}$, defining a cover $W^i$ such that \eqref{eq:nonincrease} is fulfilled replacing $\tilde\Sigma_i$ by $\tilde\Sigma_{i+1}\setminus\tilde\Sigma_i$. As a result, we have that
					\begin{equation}
						\peru(\tilde{\Sigma}_i\setminus{W^{i}};\sm{{\rm int} (\tilde \Sigma_{i+1})}) \leq \peru(\tilde{\Sigma}_i;\sm{{\rm int} (\tilde \Sigma_{i+1})}) \, \text{and} \,\leb^2(\tilde{\Sigma}_i\setminus{W^{i}}) \leq \leb^2(\tilde{\Sigma}_i)
						\label{eq:nonincrease2}
					\end{equation}

					Then, we define
					$$
					\boldsymbol{\hat\Sigma^c} = \{\hat\Sigma^c_i\}_{i=1}^{n} :=\begin{cases}
						\tilde\Sigma_i\setminus W^i & \text{if $i\in\{m_1, ... , i_0 - 1\}$}\\
						\tilde\Sigma_i\cup W^i & \text{if $i\in\{i_0, ... , m_2\}$}\\
						\Sigma_i & \text{otherwise}
					\end{cases}
					$$
					
					which satisfies that $\Gcal^n(\boldsymbol{\hat\Sigma^c}) \leq \Gcal^n(\boldsymbol{\tilde\Sigma})$ given the behaviour of $g_{c_i, c_{i+1}}$ at each index $i$ and the above inequalities \eqref{eq:nonincrease} and \eqref{eq:nonincrease2}. Moreover, their components satisfy that their boundaries are rectangular at cell $C$. Adjusting the covers $W^i$, we assure the inclusion condition $\hat\Sigma^c_i\subset\hat\Sigma^c_{i+1}$. If we repeat this process on each cell of $G_f$ adding to the changes on $\boldsymbol{\tilde\Sigma}$, we define a collection $\boldsymbol{\hat\Sigma}{\in \mathcal P_n^*(\Omega)}$ such that their components have rectangular boundaries and $$\Gcal^n(\boldsymbol{\hat{\Sigma}}) \leq \Gcal^(\boldsymbol{\tilde\Sigma}).$$
					
					\item \textit{Aligning}: All in a row, we will define a {collection}, denoted by $\boldsymbol{\Sigma}^*{\in \mathcal P_n^*(\Omega)}$, whose components are rectilinear polygons on {$\mathcal F(G_f)$}. For that, we modify the boundaries of the components of $\boldsymbol{\hat\Sigma}$ averting {an} increase of {the} energy at $\Gcal^n$.
					
					First, {for each $i=1,\ldots,n-1$}, we consider the minimal grid, denoted by $G^i$, which contains $G_f$ and $\partial\hat\Sigma_{i+1}$. Following a similar strategy to \cite[Lemma 2, Step 3]{MR3709886}, we modify $\hat\Sigma_1$ transporting segments of $\partial\hat\Sigma_{1}$ into $G^1$ while the variation of area and $\ell$--perimeter is controlled. This process, carried out segment by segment, provides a new $\hat\Sigma_{(1),1}\subset\hat\Sigma_2$ such that $\partial \hat\Sigma_{(1),1} \subset G^1$ and $\Gcal^n(\boldsymbol{\hat\Sigma}_{(1)})\leq \Gcal^n(\boldsymbol{\hat\Sigma})$, where $\boldsymbol{\tilde\Sigma}_{(1)}$ is equal to $\boldsymbol{\tilde\Sigma}$ with the exception of first element, which is replaced by $\hat\Sigma_{(1),1}$
					
					Now, we take $\hat\Sigma_2$ and we repeat the previous method but with a slight tweak. In this case, if we move a segment $s$ of $\partial\Sigma_2$ such that $s\subset\partial\hat\Sigma_{(1), 1}$, we also modify the respective part of the boundary of $\partial\hat\Sigma_{(1),1}$ in consequence. Then, this procedure gives us $ \Sigma_{(2),1}$ and $\Sigma_{(2), 2}$, variants of $\Sigma_{(1),1}$ and $\hat\Sigma_2$, such that $\partial\hat\Sigma_{(2),1}, \partial \hat\Sigma_{(2), 2}\subset G^2$ and  $\Gcal(\boldsymbol{\hat\Sigma}_{(2)}) \leq \Gcal(\boldsymbol{\hat\Sigma})$,  where $\boldsymbol{\tilde\Sigma}_{(2)}$ is equal to $\boldsymbol{\tilde\Sigma}$ with the exception of first two elements, which are replaced by $\hat\Sigma_{(2),1}$ and $\hat\Sigma_{(2),2}$. In addition, we note that the previous inclusion is fulfilled by the condition $c_i>c_{i+1}$. Rehashing this scheme for each $i$ up to $n-1$, we have defined a collection $\boldsymbol{\hat\Sigma}_{(n-1)}$ whose components have their boundaries in $G^{n-1}$.  {Moreover,  since} $\partial\hat\Sigma_n = \partial\Omega \subset G_f$, {we get} $G^{n-1} = G_f$. As a result, $\boldsymbol{\Sigma}^* := \{\Sigma_i^*\}_{i=1}^n := \boldsymbol{\hat\Sigma}_{(n-1)}{\in \mathcal P_n^*(\Omega)}$ is a collection which satisfies that
					$$
					\Gcal^n(\boldsymbol{\Sigma}^*) \leq \Gcal^n(\boldsymbol{\Sigma}) + \varepsilon \quad \text{and} \quad {\Sigma_i^*\subset \mathcal F(G_f)} \quad \text{for each $i$}.
					$$}
			\end{enumerate}
			The proof finishes as in \cite[Theorem 3]{MR3709886} and thus, we omit it.
		\end{proof}
		
		{\begin{remark}
				Observe that given $(\boldsymbol{\tilde\Sigma},{\bf c})$ a minimizer to $\mathcal G_\mu^n$, then $\boldsymbol{\tilde\Sigma}$ is a minimizer to $\mathcal{G}^n_\mu(\cdot, \boldsymbol{c})$ and that ${\bf c}$ can be reordered in such a way that $1\geq c_i\geq c_{i+1}\geq 0$ for all $i=1,\ldots,n-1$. Therefore, by Theorem \ref{PCR}, we obtain that there exists $\boldsymbol{\Sigma}\in\mathcal P^*_n(\Omega)$ such that $\boldsymbol{\Sigma}$ is a rectangular partition of $\Omega$ whose components have their boundary on $G_f$ and $(\boldsymbol{\Sigma},{\bf c})$ is also a minimizer to $\mathcal G_\mu^n$.
			\end{remark}
			
			We next show that the classical two step algorithm also leads to a PCR {function} whose components have boundaries in $G_f$, at any iteration. The algorithm reads as follows:
			
			\begin{algorithm}[H]
				\caption{$\mathcal G^n_\mu$ approximate minimizer.}
				\label{alg:CV-n_2step}
				\begin{algorithmic}
					\STATE{\textbf{Initiation:} $k = 0$, $\mu>0$, ${\bf c}=\{c_i\}_{i=1}^{n-1}\in[0,1]^{n-1}$ s.t. $c_i > c_{i+1}$.}
					\STATE{$\boldsymbol{\Sigma}^0 \gets$ $\mathcal G^n_\mu(\cdot,{\bf c})$ PCR minimizer via Theorem \ref{PCR}.}
					\vspace{0.2cm}
					
					\WHILE{$\sum_{i=0}^n |\Sigma_i^k \bigtriangleup \Sigma_{i}^{k-1}|^2 > \varepsilon_{tol} \, \land \, k < n_{max}$}
					\vspace{0.2cm}
					\STATE{$c_i \gets \frac{1}{\left|\Sigma^k_{i}\setminus\Sigma^k_{i-1}\right|} \int_{\Sigma^k_{i}\setminus\Sigma^k_{i-1}} f \, dx$,\quad 	$k \gets k + 1$;}
					\STATE{$\boldsymbol{\Sigma}^k\gets$ $\mathcal G^n_\mu(\cdot,{\bf c})$  PCR minimizer via Theorem \ref{PCR};}
					
					\vspace{0.2cm}
					\ENDWHILE
					\RETURN{$(\boldsymbol{\Sigma}^k,{\bf c})$}
					
				\end{algorithmic}
			\end{algorithm}
			
			\bigskip	
			In order to be able to apply Theorem \ref{PCR}, we only need to show next result.}
		\begin{proposition}\label{prop3}
			Let ${\bf c}=\{c_i\}_{i=1}^n$ satisfy $1\geq c_i\geq c_{i+1}\geq 0$, and let $\boldsymbol{\Sigma}$ be the $\pcr(f)$ minimizer of $\mathcal{G}^n_\mu(\cdot,{\bf c})$ obtained in Theorem \ref{PCR}. Then, $\tilde c_{i}> \tilde{c}_{i+1}$, for $i=1,\ldots,n-1$ with $\tilde c_i$ defined as $$\tilde c_i:=\frac{1}{|\Sigma_{i}\setminus\Sigma_{i-1}|}\int_{\Sigma_i\setminus\Sigma_{i-1}}f(x)\,dx.$$
		\end{proposition}
		
		\begin{proof}
			First of all, note that, given $\Sigma_{i-1}$ and $\Sigma_{i+1}$, $\Sigma_i$ is a minimizer of $$\peru(\Sigma_{i-1},\sm{{\rm int}}(E))+\peru(E,\sm{{\rm int}}(\Sigma_{i+1}))+\mu (c_i-c_{i+1})\int_{E}(c_{i}+c_{i+1}-2f)\,dx, $$ for all sets $\Sigma_{i-1}\subset E\subset \Sigma_{i+1}$. Note that the minimizer does not depend on the behaviour {of $f$} on $\Sigma_{i-1}$ and on $\Sigma_{i+1}$. Therefore, $\Sigma_i\setminus\Sigma_{i-1}$ minimizes in $\Sigma_{i+1}\setminus\Sigma_{i-1}$, the following functional:
			\begin{equation}\label{problemini}
				\peru(E,\sm{{\rm int}}(\Sigma_{i+1}\setminus\Sigma_{i-1}))+\mu(c_i-c_{i+1})\int_E (c_i+c_{i+1}-2f)\,dx.
			\end{equation}
			Since $\emptyset$ is admissible, and $\Sigma_{i}\setminus\Sigma_{i-1}$ is a minimizer, it follows that$$\int_{\Sigma_{i}\setminus\Sigma_{i-1} } (c_i+c_{i+1}-2f)\,dx{\leq}0\,.$${In case of equality, the perimeter is $0$ and then, $\Sigma_i$ coincides either with $\Sigma_{i-1}$ or $\Sigma_{i+1}$) and that phase will be removed passing from $n-$ to $(n-1)-$phases. Therefore, we can suppose that} $$\tilde c_i > \frac{c_i+c_{i+1}}{2}$$
			Therefore, if $\tilde c_{i}\leq \tilde c_{i+1}$, it follows that $$\int_{\Sigma_{i+1}\setminus\Sigma_i}(c_i+c_{{i+1}}-2f)\,dx< 0.$$This implies that $\Sigma_{i+1}{\setminus \Sigma_{i-1}}$ has strictly less energy in \eqref{problemini} than $\Sigma_i{\setminus \Sigma_{i-1}}$, a contradiction.
		\end{proof}

		\section{Relationship between ACV and Truncated AROF models }\label{sec:trof}
		
		{This section is aimed at showing the relation between $\cv_{\varphi,\boldsymbol\mu}^n$ and ${\rm TROF}_{\varphi, \lambda}^n$ minimizers with respect to $|\cdot|_\varphi$, functionals defined in \eqref{eq:cvphi} and \eqref{eq:trofphi}, respectively. Throughout this section, ${\rm TROF}_\lambda^n$ will denote the anisotropic version of ${\rm TROF}_{\varphi, \lambda}^n$.
			
			In the anisotropic 2-phases case, if $\boldsymbol\mu$ is a constant vector, we have that the minimizer of ${{\rm TROF}_{\lambda}^2(\{\emptyset, \ \cdot \ ,\Omega\}, \tau)}$ is provided by the upper level set $\Sigma_{\tau}:=\{x\in\Omega : w_\lambda(x)>{\tau}\}$ where $w_\lambda$ is the solution of $AROF_\lambda$, by analogy on the result \cite[Proposition 2.6]{CHAMBOLLE-CASELLES}. In addition, by Theorem \ref{thm:CV1sol}, we have that
			
			\begin{corollary}\label{col:equiv2phases}
				{Let  $w$ be the minimizer of $\rofa_{\lambda}$ and $\Sigma_\tau$ the $w$'s level set solution of ${\rm TROF}_{\lambda}^2(\{\emptyset, \ \cdot \ ,\Omega\}, \tau)$ such that $\tau$ is defined as
					\begin{align}
						\tau = \frac{c_1 + c_2}{2}, \, \, \text{with} \, \, c_1 := \frac{\int_{\Sigma_\tau} f\, dx}{|\Sigma_\tau|}, \, \, c_2 := \frac{\int_{\Omega \setminus\Sigma_\tau} f\, dx}{|\Omega\setminus\Sigma_\tau|}.
					\end{align}
					{ Then, $(\Sigma_\tau, c_1, c_2)$ is a minimizer of $\acv_\mu$ if $\lambda = {\mu}(c_1^* - c_2^*)$ and $(\chi_{\Sigma_\tau}, c_1, c_2)$ is a minimizer of $\mathcal G_\mu$.}}
			\end{corollary}
			
			In view of this result, which connects $AROF_\lambda, ACV_\mu$ and the anisotropic ${\rm TROF}_{{1},\lambda}^2$, it is relevant to propose similar relationships on {the} multiphase case. However, {as we will show,} the relationship between $CV_{\varphi, \boldsymbol{\mu}}^n$ and ${\rm TROF}_{\varphi,\lambda}^n$, presented in \cite{Cai} in the isotropic case, does not {hold, in general}.} The connection is presented with respect to \sm{a generic anisotropy} $|\cdot|_\varphi$ next:\\
		
		\textbf{Relationship} \textbf{CV}$_{\boldsymbol\varphi,\boldsymbol\mu}^{\boldsymbol n}$ \textbf{- TROF}$_{\boldsymbol\varphi,\boldsymbol\lambda}^{\boldsymbol n}$:
		{\it Let $(\boldsymbol\Sigma, \boldsymbol{\tau})$ be a pair in $\mathcal P_{n}^*(\Omega)\times[0,1]^n$ such that
			\begin{equation*}
				{\rm TROF}_{\varphi, \lambda}^{n}(\boldsymbol\Sigma, \boldsymbol \tau) \leq {\rm TROF}_{\varphi, \lambda}^{n}(\boldsymbol\Lambda, \boldsymbol \tau)
			\end{equation*}
			for any feasible $\boldsymbol\Lambda \in \mathcal P_{n}^*(\Omega)$; and
			\begin{equation*}\tau_i = \frac{c_{i+1} + c_i}{2} \ \ \text{with} \ \ c_i := \int_{\Omega_i} \frac{f}{|\Omega_i|}\, dx, \, \, \Omega_i := \Sigma_{i}\setminus\Sigma_{i-1},
			\end{equation*}
			for each $i\in\{1, ..., n-1\}$. Provided that $c_i>c_{i+1}$, we define $\boldsymbol\mu := \{\mu_i\}_{i=1}^n$ as follows
			\begin{equation}\label{def:mu}
				\begin{split}
					&\mu_1 = \frac{\lambda}{2(c_1 - c_2)}, \, \,  \mu_n = \frac{\lambda}{2(c_{n-1} - c_{n})},\\
					&\mu_i =\frac{\lambda(c_{i-1}-c_{i+1})}{2(c_{i-1} - c_i)(c_i - c_{i+1})} \quad \text{for $1<i<n$},
				\end{split}
			\end{equation}
			{The possible relationship is this one: Letting } $\boldsymbol\Omega :=\{\Omega_i\}_{i=1}^{n}$ and $\boldsymbol c := \{c_i\}_{i=1}^{n}$, {is it true that }
			\begin{equation}\label{eq:jointquasiminimiser}
				\cv_{\varphi,\boldsymbol\mu}^n(\boldsymbol\Omega, \boldsymbol c) \in \argmin_{{\boldsymbol \Lambda\in\mathcal P_n(\Omega) }}  \cv_{\varphi,\boldsymbol\mu}^n(\boldsymbol\Lambda, \boldsymbol c) \, ?
			\end{equation}
		}
		To show {that this is not the case}, we need some preliminary results. These results are an adaptation of the results in \cite[Section 9]{BelletiniCasellesNovaga}. Since our aim is to show some very particular examples, we only {consider} some specific settings.
		
		\begin{definition}
			Let $A\subset\mathbb{R}^2$ be a bounded open set with Lipschitz boundary. We say that $A$ is calibrable with respect to the anisotropy $\varphi$ if there exists $\z\in { X_A}$ such that $|\z|^*_\varphi\leq 1$ a.e. in $A$,  ${\rm div} \z$ is constant in $A$ and $[\z,\nu^A]=-|\nu^A|_\varphi^*$ at $\partial A$.
		\end{definition}
		
		\begin{lemma}\label{lem:2}
			Let $A\subset\mathbb{R}^2$ be a bounded open set with Lipschitz boundary, and let $\Omega:= {B_{\varphi^*}}(0;R)$ be such that ${\overline A\subset\Omega}$. Then, there exists $\z\in X_{\Omega\setminus \overline A}$ such that $|\z|_\varphi^*\leq 1$ a.e. in $\Omega\setminus \overline A$, ${\rm div} \,\z$ is constant in $\Omega\setminus \overline A$ and $[\z,\nu^{\Omega}]=0$ at $\partial \Omega$ and $[\z,\nu^A]=-|\nu^A|_\varphi^*$ at $\partial A$ if and only if
			\begin{equation}\label{conditionB_R}
				\Omega\setminus \overline A \in \argmin_{E\in \mathcal E}\left\{\frac{\per_\varphi(E; \ \Omega\setminus \overline A)-\mathcal H^1(\partial^*E\cap \partial A)}{|E|}\right\},
			\end{equation}
			with $\mathcal E:=\{ E\subset \Omega\setminus \overline A \ : \ \per_\varphi(E)<+\infty\,, |E|>0\}$.
		\end{lemma}
		\begin{proof}
			The proof follows the same ideas as the proof of \cite[Theorem 5]{BelletiniCasellesNovaga}.
			
			Suppose first that there is a vector field $\z$ satisfying the hypothesis of the Lemma. Then, by Grenn-Gauss theorem,
			$$
			\int_{\Omega \setminus\overline A}{\rm div \ }\z\, dx =-\int_{\partial(\Omega\setminus\overline A)}[\z,\nu]\,d\mathcal H^1=\per_\varphi(A).
			$$
			
			Therefore, ${\rm div \ }\z=\frac{\per_\varphi(A)}{|\Omega\setminus\overline A|}$. Applying once again Green-Gauss theorem, for any {$E\in \mathcal{E}$}, we obtain
			$$
			|E|\frac{\per_\varphi(A)}{|\Omega \setminus\overline A|}=\int_{E}{\rm div \ }\z\,dx =-\int_{\partial E}[\z,\nu]\,d\mathcal H^1\geq \mathcal H^1(\partial^*E\cap \partial A)-\per_\varphi(E; \ \Omega\setminus\overline A),
			$$
			which shows the first implication.
			
			\medskip
			Let us suppose now that $\Omega\setminus\overline A$ is a minimizer of the functional
			$$
			E\mapsto \frac{\per_\varphi(E;\ \Omega\setminus \overline A)-\mathcal H^1(\partial^*E\cap \partial A)}{|E|},$$among all sets in $\mathcal E$. We now define the functional
			$$
			F(\boldsymbol{\xi}):=\int_{\Omega\setminus A}({\rm div}\boldsymbol{\xi})^2\,dx, \quad \text{with} \, \, \boldsymbol{\xi} \in X_{\Omega\setminus \overline A}.
			$$
			Therefore, arguing as in \cite[Proposition 6.1 and Theorem 6.7]{bellettininovagapaolini1} and \cite[Proposition 3.5 and Theorem 5.3]{bellettininovagapaolini2}, one can prove that the following variational problem has a solution with unique divergence:
			\[
			\min\left\{F(\boldsymbol{\xi}) \, : \, \begin{array}{l}
				\boldsymbol{\xi}\in X_{\Omega\setminus \overline A}\,, \quad |\boldsymbol{\xi}|_\varphi^*\leq 1 {\rm \ a.e. \ in \ } \Omega \setminus \overline A\\
				{ [\boldsymbol{\xi},\nu^{\Omega}]=0 {\rm \ at \ } \partial \Omega {\rm\quad  and \quad } [\boldsymbol{\xi},\nu^A]=-|\nu^A|_\varphi^* {\rm \ at \ }\partial A }
			\end{array}
			\right\}\,.
			\]
			
			Moreover, given any minimizer $\boldsymbol{\xi}_{\min}$, ${\rm div }\boldsymbol{\xi}_{\min}$ $\in L^\infty(\Omega\setminus\overline A)\cap \bv(\Omega\setminus\overline A)$  and, letting $Q_\mu$ be the $\mu$ upper level set of ${\rm div }\boldsymbol{\xi}_{\min}$ in $\Omega\setminus\overline A$, $Q_\mu$ has finite perimeter, and
			\begin{equation}\label{Q_mu}
				\int_{Q_\mu}{\rm div \ }\boldsymbol{\xi}_{\min} \,dx=\mathcal H^1(\partial^*Q_\mu\cap\partial A)-\per_\varphi(Q_\mu; \ \Omega\setminus\overline A).
			\end{equation}
			
			Were ${\rm div \, }\boldsymbol{\xi}_{\min}$ not constant in $\Omega\setminus\overline A$ and equal to $\frac{\per_\varphi(A)}{|\Omega\setminus\overline A|}$, there will be $\mu_0>\frac{\per_\varphi(A)}{|\Omega\setminus\overline A|}$ such that $Q_{\mu_0}$ is nonempty. Therefore, by \eqref{Q_mu}, we obtain that $\Omega\setminus\overline A$ cannot be a minimizer of the above functional.
		\end{proof}
		
		\begin{theorem}\label{lem:3}
			Let $A\subset\mathbb{R}^2$ be a bounded open set with Lipschitz boundary such that $A=\cup_{i=1}^m C_i$ with $C_i$ disjoint, convex and calibrable, and let $\Omega := {B_{\varphi^*}}(0;R)$ be such that ${A\subset\subset \Omega}$. Let $1\leq k\leq m$, $\{i_1,\ldots,i_k\} \subset \{1,\ldots,m\}$ and consider the following variational problem
			$$
			(P)_{i_1,\ldots,i_k}:=\min_{E\in \mathcal{E}_{i_1,\ldots,i_k}} \left\{\per_\varphi(E; \ \Omega)+\frac{\per_\varphi(A)}{\left|\Omega\setminus\overline A \right|}|E\setminus\overline A| \right\}\,,
			$$
			with $\mathcal{E}_{i_1,\ldots,i_k}:=\left\{\cup_{j=1}^k C_{i_j}\subseteq E\subseteq \Omega\setminus \cup_{j=k+1}^m C_{i_j} \ : \ \per_\varphi(E)<+\infty \right\}$.
			
			Then, there exists $\z\in X_{\Omega\setminus \overline A}$ such that $|\z|_\varphi^*\leq 1$ a.e. in $\Omega\setminus \overline A$, ${\rm div}\z$ is constant in $\Omega\setminus \overline A$ and $[\z,\nu^{\Omega}]=0$ at $\partial \Omega$ and $[\z,\nu^A]=-|\nu^A|_\varphi^*$ at $\partial A $
			if and only if \,  $\bigcup_{j=1}^k C_{i_j}$ is a solution to $(P)_{i_1,\ldots,i_k}$ for any $\{i_1,\ldots,i_k\} \subset \{1,\ldots,m\}$.
		\end{theorem}
		
		\begin{proof}
			
			Suppose first that there exist a vector field $\z$ as in the hypothesis and let us consider {$E\in \mathcal{E}_{i_1,\ldots,i_k}$} and $D:=E\setminus \overline A$.
			Then, by Green-Gauss theorem,
			$$
			\frac{\per_\varphi(A)}{\left|\Omega\setminus\overline A\right|}|D|=\int_{D}{\rm div \ }\z\,dx\geq -\per_\varphi(E;\ \Omega\setminus\overline A)+\mathcal H^1(\partial^* D\cap\bigcup_{j=1}^k\partial C_{i_j}).
			$$
			Therefore, $$
			\sum_{j=1}^k\per_\varphi(C_{i_j})\leq \sum_{j=1}^k\per_\varphi(C_{i_j}) +\frac{\per_\varphi(A)}{\left|\Omega\setminus\overline A\right|}|D|+ \per_\varphi(E; \Omega\setminus\overline A)$$$${ -} \mathcal H^1(\partial^* D\cap\bigcup_{j=1}^k\partial C_{i_j})\leq \frac{\per_\varphi(A)}{\left|\Omega\setminus\overline A\right|}|D|+\per_\varphi(E;\ \Omega)\,,
			$$
			thus showing that $\bigcup_{j=1}^k C_{i_j}$ is a solution to $(P)_{i_1,\ldots,i_k}$. 
			
			Now, let us suppose that $\bigcup_{j=1}^k C_{i_j}$ is a solution to $(P)_{i_1,\ldots,i_k}$ for any $\{i_1,\ldots,i_k\} \subset \{1,\ldots,m\}$.   We only need to show that, in this specific case of $A$, \eqref{conditionB_R} holds.
			Let $E$ be a set of finite perimeter such that $E\subseteq \Omega\setminus \overline A$ and let $C_{i_j}$ for $j=1,\ldots,k$ be the ones such that $\partial E\cap C_{i_j}\neq \emptyset$. Then, by minimality,
			$$
			\sum_{j=1}^k\per_\varphi(C_{i_j})\leq \per_\varphi\left(E\cup\bigcup_{j=1}^k C_{i_j}; \ \Omega\right)+\frac{\per_\varphi(A)}{\left|\Omega\setminus\overline A \right|}|E|
			$$
			$$
			\leq \per_\varphi(E; \ \Omega)+\sum_{j=1}^k\per_\varphi(C_{i_j})-\mathcal H^1(\partial^*E\cap\partial A)+\frac{\per_\varphi(A)}{\left|\Omega\setminus\overline  A \right|}|E|.
			$$
			Therefore,
			$$
			\per_\varphi(E; \ \Omega)-\mathcal H^1(\partial^*E\cap\partial A)+\frac{\per_\varphi(A)}{\left|\Omega\setminus\overline A \right|}|E|\geq 0.
			$$
			This implies \eqref{conditionB_R}.
		\end{proof}
		
		\begin{corollary}
			\label{col:1}
			Let $A\subset\mathbb{R}^2$ be a bounded open set with Lipschitz boundary such that $A=\cup_{i=1}^m C_i$ with $C_i$ disjoint, convex and calibrable, and let $\Omega := B_{\varphi^*}(0;R)$ be such that ${\overline A\subset\Omega}$. Suppose that
			\begin{equation}\label{eq:distcondition}
				{\rm dist }_\varphi\left(C_i,\bigcup_{j\neq i}C_j\cup  {\left( \mathbb{R}^2\setminus\overline\Omega \right)}\right)>\per_\varphi(C_i)\,\quad {\rm for \ any \ }i=1,\ldots,m.
			\end{equation}
			Let $f:=\sum_{i=1}^m \alpha_i\chi_{C_i}$ with $\alpha_i>0$. Then, for
			\begin{equation}\label{eq:lambdacondition}
				\lambda\geq \max_{i=1,\ldots,m} \left\{ \frac{1}{\alpha_i}\left(\frac{\per_\varphi(A)}{|\Omega\setminus\overline A|}+\frac{\per_\varphi(C_i)}{|C_i|}\right)\right\},
			\end{equation}
			the solution to {{\rm ROF}$_{\varphi,\lambda}$} model is given by $$u=\sum_{i=1}^m\left(\alpha_i-\frac{\per_\varphi(C_i)}{\lambda|C_i|}\right)\chi_{C_i}+\frac{\per_\varphi(A)}{\lambda|\Omega\setminus\overline A|}\chi_{\Omega\setminus A}.$$
		\end{corollary}
		\begin{proof}
			Since each $C_i$ is calibrable, we can construct a vector field $\z_i\in X_{C_i}$ such that ${\rm div \ }\z_i=\frac{\per_\varphi(C_i)}{|C_i|}$ and $[\z_i,\nu^{C_i}]=-|\nu^{C_i}|_\varphi^*$. Moreover, it is easy to show that the hypothesis in Theorem \ref{lem:3} are satisfied. Therefore, we can construct a vector field $\z_{\rm out}\in X_{\Omega\setminus\overline A}$ such that
			${\rm div \ }\z_{\rm out}=\frac{\per_\varphi(A)}{|\Omega\setminus\overline A|}$, $[\z_{\rm out},\nu^{C_i}]=-|\nu^{C_i}|_\varphi^*$ and $[\z_{\rm out},\nu^{\Omega}]=0.$
			{Furthermore}, observe that, in this case $\rest{u}{C_i}\geq \rest{u}{\Omega\setminus\overline A}$. Therefore, it is easy to show that, considering
			$$
			\z:=\sum_{i=1}^m \z_i\chi_{C_i}+\z_{\rm out}\chi_{\Omega\setminus\overline A},
			$$
			then $-{\rm div \ }\z\in\partial |Du|_\varphi$ and therefore, $u$ is the solution of the {ROF$_{\varphi,\lambda}$} model {by \eqref{eq:ELphiROF}.}
		\end{proof}

		\begin{example}
			{{In this example, using a case of $3$-phases segmentation, we show that the relationship \eqref{eq:jointquasiminimiser} does not hold, in general.}} Let $\Omega:= B_{\varphi^*}((0,0), R)$ and $f = \chi_{C_1} +  \chi_{C_2} + \frac{1}{2}\chi_{C_3}$ {with}  $C_1 := B_{\varphi^*}((0, -L),1)$, $C_2 := B_{\varphi^*}((0,0), (1/2))$ and $C_3 := B_{\varphi^*}((0,L), 2)$. We {take} $R$ and $L$ large enough so that $\lambda = 10$ is feasible with respect to \eqref{eq:lambdacondition}, $A:=\cup_{i=1}^3 C_i$ satisfies that $\overline{A} \subset \Omega$  and  $C_i$ fulfils \eqref{eq:distcondition} for each $i$. Then, by Corollary \ref{col:1}, the solution of ${\rm ROF}_{\varphi, \lambda}$ problem is exactly given by
			\begin{equation}\label{eq:phiROFsol}
				w=\frac{8}{10}\chi_{C_1}+\frac{6}{10}\chi_{C_2}+\frac{4}{10}\chi_{C_3} + \frac{\per_\varphi(A)}{10\left|\Omega\setminus \overline{A}\right|}\chi_{\Omega\setminus \overline{A}}
			\end{equation}
			
			Now, we define $c_1, c_2$ and $c_3$ as follows:
			\begin{align*}
				c_1 = \int_{C_1} \frac{f}{|C_1|}\, dx \,, \quad c_2 = \int_{C_2 \cup C_3} \frac{f}{|C_2 \cup C_3|}\, dx \,, \quad c_3 = \int_{\Omega\setminus \overline A} \frac{f}{\left |\Omega\setminus \overline A \right|}\,dx \, .
			\end{align*}
			where $c_1 = 1, c_2 = 9/17$ and $c_3 = 0$. By enlarging $R$ if necessary, we \sm{assume} that the latter term of \eqref{eq:phiROFsol} is bounded on top by $c_2/2$. Thus, we may assure that
			\begin{align}
				\Sigma_1 &:= \left\{x\in\Omega : w >\frac{c_1 + c_2}{2}\right\}=C_1\,,\quad
				\Sigma_2 := \left\{x\in\Omega : w >\frac{c_2 + c_3}{2} \right\}=A \,.
			\end{align}
			
			If we suppose {that \eqref{eq:jointquasiminimiser} is true, then} $\{\Omega_i\}_{i=1}^3:=\{\Sigma_1, \, \Sigma_2\setminus\Sigma_1,\, \Omega\setminus\Sigma_2\} = \{C_1, \, C_2\cup C_3, \, \Omega\setminus \overline{A}\}$ satisfies that
			\begin{equation}\label{eq:quasiminimizersvcv}
				\cv_{\varphi,\boldsymbol{\mu}}^3(\{\Omega_i\}_{i=1}^3, \{c_i\}_{i=1}^3) \leq \cv_{\varphi,\boldsymbol{\mu}}^3(\{\Lambda_i\}_{i=1}^3, \{c_i\}_{i=1}^3)\,, \qquad \forall\{\Lambda_i\}_{i=1}^3\in\mathcal{P}_3(\Omega),
			\end{equation}
			where $\boldsymbol{\mu}= \{85/8,\, 1445/72, \, 85/9\}$,  as in \eqref{def:mu}. However, if we define the following disjoint partition $\{\Omega_i^*\}_{i=1}^3 := \{C_1\cup C_2, \, C_3, \, \Omega\setminus \overline{A}\}$, we have that
			$$
			\cv_{\varphi, \boldsymbol{\mu}}^3(\{\Omega^*_i\}_{i=1}^3, \{c_i\}_{i=1}^3) = \frac{725}{72}|C_1| < \cv_{\varphi, \boldsymbol{\mu}}^3(\{\Omega_i\}_{i=1}^3, \{c_i\}_{i=1}^3) = \frac{733}{72}|C_1|,
			$$ which leads us to a contradiction with respect to \eqref{eq:quasiminimizersvcv}. It should be noted that this reasoning is valid for any anisotropy $|\cdot|_\varphi$.
		\end{example}

	{Finally, we would like to point out a relationship between $\rofa_\lambda$ and ${\rm TROF}_{ \lambda}^n$. First we note that, for a fixed $\boldsymbol\tau\in[0,1]^{n-1}$, any $\boldsymbol{\Sigma}\in\mathcal P_{n}^*(\Omega)$ satisfies this expression
		\begin{equation}
			{\rm TROF}_{\varphi,\lambda}^n(\boldsymbol\Sigma, \boldsymbol\tau) = \sum_{i=1}^{n-1} {\rm TROF}_{\varphi, \lambda}^2(\{\emptyset, \Sigma_i, \Omega\}, \tau_i)\,.
		\end{equation}
		Thus, Corollary \ref{col:equiv2phases} implies the following result.
		
		\begin{corollary}\label{col:algorithm3}
			Let $w$ be the $\rofa_\lambda$ minimizer. Then, for a $\boldsymbol{\tau}\in[0,1]^{n-1}$ such that $\tau_{i}>\tau_{i+1}$, the minimizer $\boldsymbol{\Sigma}\in\mathcal P_{n}^*(\Omega)$ of ${\rm TROF}_{\lambda}^n(\cdot, \boldsymbol\tau)$ is provided by
			\begin{equation}
				\Sigma_i = \{x\in\Omega : w(x)>\tau_i\} \quad \text{for} \,\, i=1, ..., n-1.
			\end{equation}
		\end{corollary}
		
		This fact allows us to define an algorithm that, by calculating the minimizer of $\rofa_\lambda$, {we can obtain} a minimizer of the functional ${\rm TROF}_{\lambda}^n$. In Section \ref{sec:applications} we show that it can be an advantageous segmentation tool in certain situations. The algorithm, defined as the previous ones, is the following one:

		\begin{algorithm}[H]
			\caption{${\rm TROF}^n_{\lambda}$ approximate minimizer}
			\label{alg:TAROFnsketch}
			\begin{algorithmic}
				
				\STATE{\textbf{Initiation:} $k = 0$, $\{\tau_i\}_{i=1}^{n-1}\in[0,1]^{n-1}$ s.t. $\tau_i > \tau_{i+1}$.}
				\STATE{$w_0 \gets$ $\rofa_\lambda$ minimizer.}
				\STATE{$ \boldsymbol\Sigma^0 \gets$ ${\rm TROF}^n_{\lambda}$ minimizer via Corollary \ref{col:algorithm3}.}
				\vspace{0.2cm}
				\WHILE{$\sum_{i=0}^n |\Sigma_i^k \bigtriangleup \Sigma_{i}^{k-1}|^2 > \varepsilon_{tol} \, \land \, k < n_{max}$}
				\vspace{0.2cm}
				\STATE{$\tau_i \gets \frac{1}{\left|\Sigma^k_{i+1}\setminus\Sigma^k_{i}\right|} \int_{\Sigma^k_{i+1}\setminus\Sigma^k_{i}} f \, dx$,\quad 	$k \gets k + 1$.}\\
				\vspace{0.1cm}
				\STATE{$\boldsymbol\Sigma^k \gets$ ${\rm TROF}^n_{\lambda}$ minimizer via Corollary \ref{col:algorithm3}.}
				\vspace{0.2cm}
				\ENDWHILE
				\RETURN{$\boldsymbol\Sigma^k$}
			\end{algorithmic}
		\end{algorithm}

		\begin{remark}
			We should note that, although the relationship between $\cv_{\varphi,\boldsymbol\mu}^n$ and ${\rm TROF}_{\varphi,\lambda}^n$ \eqref{eq:jointquasiminimiser} does not hold, in general, the segmentation provided by the ${\rm TROF}_{\varphi,\lambda}^n$ functional provides satisfactory \sm{results} both in the isotropic case and in our case: the anisotropic one {\vp (e.g. see figure \ref{fig:example-cub})}. Moreover, in \cite{Cai} it was shown that algorithm \ref{alg:TAROFnsketch}, in the isotropic formulation, converges. {Their proof can be easily adapted in the anisotropic framework and thus we can prove that the above algorithm (directly) and Algorithm \ref{alg:ACVsketch} (applying Corollary \ref{col:equiv2phases}) converge.} Therefore, in the next section we will make use of algorithm \ref{alg:TAROFnsketch} when studying some multiphase segmentations.	
		\end{remark}

	}

	\section{Applications}\label{sec:applications}
\begin{figure}[H]
		\centering
		\begin{subfigure}[b]{0.45\textwidth}
			\centering
			\includegraphics[width=0.7\textwidth]{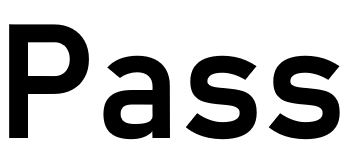}
			\caption{Original image}
			\label{fig:pass}
		\end{subfigure}
		\begin{subfigure}[b]{0.45\textwidth}
			\centering
			\includegraphics[width=0.7\textwidth]{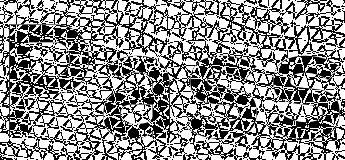}
			\caption{Corrupted image}
			\label{fig:passlinies}
		\end{subfigure}
		\begin{subfigure}[b]{0.45\textwidth}
			\centering
			\includegraphics[width=0.7\textwidth]{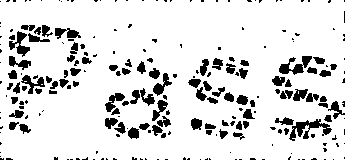}
			\caption{\cite{Cai} approach with $\mu = 2.9$}
			\label{fig:pascai}
		\end{subfigure}
		\begin{subfigure}[b]{0.45\textwidth}
			\centering
			\includegraphics[width=0.7\textwidth]{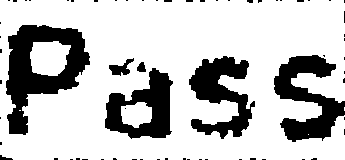}
			\caption{Alg. \ref{alg:ACVsketch} approach with $\mu = 4.3$}
			\label{fig:pasmp2}
		\end{subfigure}
		
		\caption{Comparison of $2$-phases segmentation.}
		\label{fig:example-pass}
	\end{figure}

\begin{figure}[H]
		\centering
		\begin{subfigure}[b]{0.49\textwidth}
			\centering
			\includegraphics[width=0.7\textwidth]{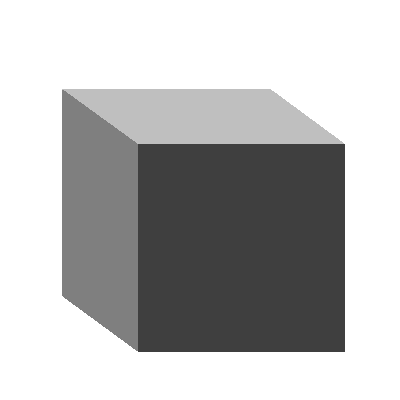}
			\caption{Original image}
			\label{fig:cub}
		\end{subfigure}
		\begin{subfigure}[b]{0.49\textwidth}
			\centering
			\includegraphics[width=0.7\textwidth]{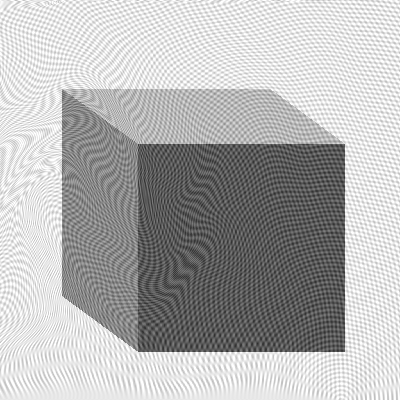}
			\caption{Corrupted image}
			\label{fig:cubmoire}
		\end{subfigure}
		\begin{subfigure}[b]{0.49\textwidth}
			\centering
			\includegraphics[width=0.7\textwidth]{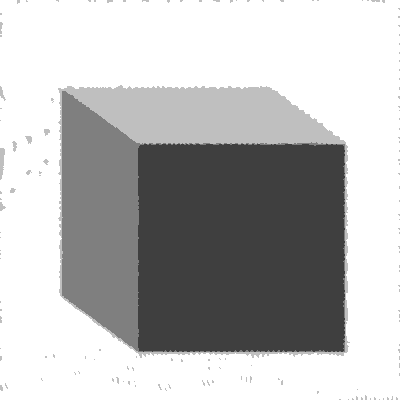}
			\caption{\cite{Cai} approach with $\lambda = 4.2$}
			\label{fig:cubcai}
		\end{subfigure}
		\begin{subfigure}[b]{0.49\textwidth}
			\centering
			\includegraphics[width=0.7\textwidth]{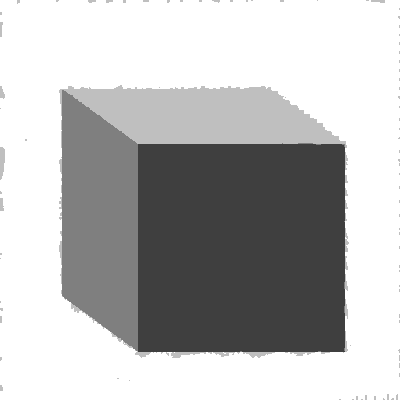}
			\caption{Alg. \ref{alg:TAROFnsketch} approach with $\lambda = 11.5$}
			\label{fig:cubmp2}
		\end{subfigure}
		
		\caption{Comparison of $4$-phases segmentations.}
		\label{fig:example-cub}
	\end{figure}

	In the previous sections,  we have \sm{looked} into an analytical approach to $\rofa_\lambda$, $\acv_\mu$, $\cv_{\varphi,\boldsymbol\mu}^n$ and ${\rm TROF}_{\varphi,\lambda}^n$, which gives us a robust tool to find solutions to these variational models. In the literature, CV and ROF models have been studied from a approximative perspective, where successful methods of resolution have been proposed which lead to an approximation of the theoretical minimizer {(see e.g. \cite{FIGUEIREDO, PARIKH, CHAMBOLLE-POCK, HE-OSHER, KEEGAN}). In this section we will show that the theoretical study of these models provides efficient segmentation algorithms and characterizations of minimizers that are of great interest. Therefore, we show the performance of algorithms \ref{alg:ACVsketch} and \ref{alg:TAROFnsketch}, comparing them with other ones based on level sets; and we present a possible application of theorems \ref{thm:CV1sol} and \ref{PCR}.\\
		
		To exemplify the performance of algorithms \ref{alg:ACVsketch} and \ref{alg:TAROFnsketch}, we compare these algorithms with the similar isotropic ones proposed in \cite{Cai}. Specifically, in figures \ref{fig:example-pass} and \ref{fig:example-cub}  we compare algorithms \ref{alg:ACVsketch} and \ref{alg:TAROFnsketch} with their isotropic counterparts.} For this, we set $n_{max}=200$, $\varepsilon_{tol}=10^{-3}$, we initialise $c_i$ and $\tau_i$ using FCM processes (see \cite{FCM}) and we approximate the $\rofa_\lambda$ and ${\rm ROF}_\lambda$ minimizers via standard Split-Bregman algorithms. Moreover, the selected examples are those segmentations that differ the least, in quadratic terms, from the original images.\\
	
	In figures \ref{fig:example-pass} and \ref{fig:example-cub}, we observe that our anisotropic model is less influenced by the \sm{applied} noise 
	than the isotropic one. Thus, it provides a better reconstruction of word \textit{Pass} (compare figures \ref{fig:pascai} and \ref{fig:pasmp2}) or of the edges of a shaded cube (compare figures \ref{fig:cubcai} and \ref{fig:cubmp2}). This fact may be correlated with a worse performance on the isotropic $TV$ on these noises. As the compared algorithms are, to some extent, analogous to ours, these examples show how the models studied in this work provide \sm{some} advantages in certain situations where isotropic versions of these do not \sm{give} adequate segmentations.\\
	
\begin{figure}[H]
		\centering
		\begin{subfigure}[b]{0.49\textwidth}
			\centering
			\includegraphics[width=0.75\textwidth]{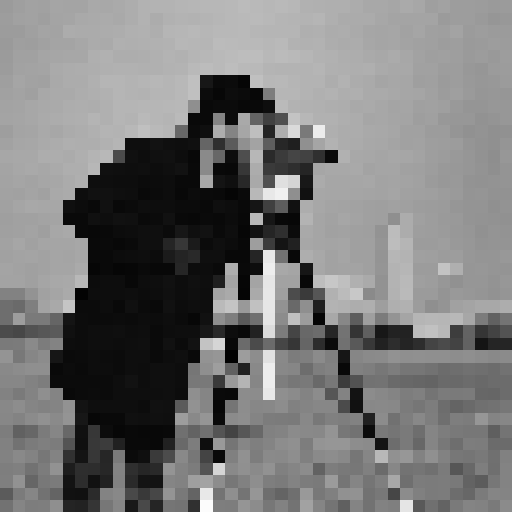}
			\caption{Original image}
			\label{fig:cman}
		\end{subfigure}
		\hfill
		\begin{subfigure}[b]{0.49\textwidth}
			\centering
			\includegraphics[width=0.75\textwidth]{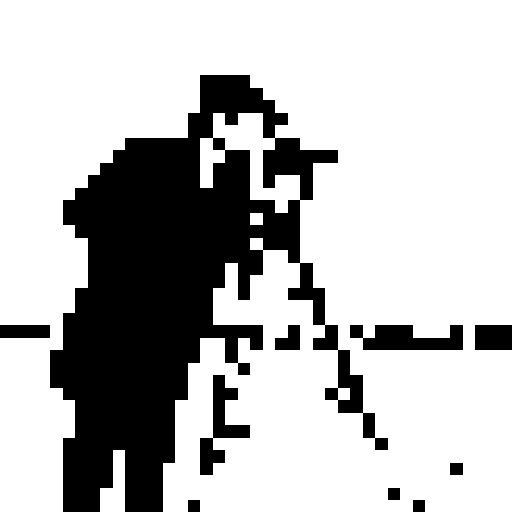}
			\caption{Exact $ACV_\mu$ minimizer}
			\label{fig:cman-orig}
		\end{subfigure}
		\hfill
		\begin{subfigure}[b]{0.49\textwidth}
			\centering
			\includegraphics[width=0.75\textwidth]{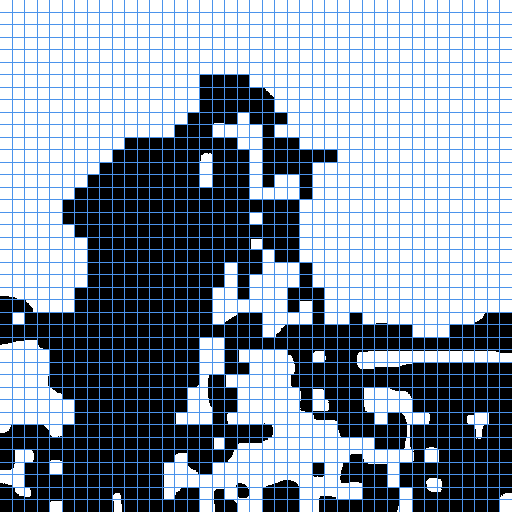}
			\caption{Anisotropic \cite{Getreur} approach}
			\label{fig:cman-getreur}
		\end{subfigure}
		\begin{subfigure}[b]{0.49\textwidth}
			\centering
			\includegraphics[width=0.75\textwidth]{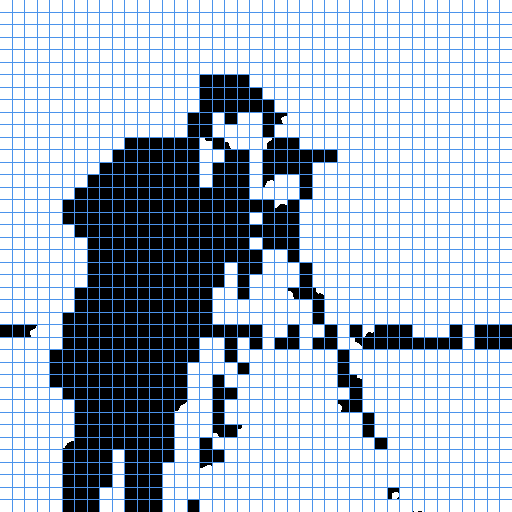}
			\caption{Algorithm \ref{alg:ACVsketch}'s approach}
			\label{fig:cman-mp2}
		\end{subfigure}
		\caption{ Comparison between two $ACV_\mu$ approaches, for $\mu = 10$. Segmentations and input image grid are in black-white and blue, respectively.}
		\label{fig:example-cman}
	\end{figure}

	Finally,  we discuss the application of Theorem \ref{thm:CV1sol} and \ref{PCR} as a reliability analysis tool. According to the previous results, we could consider that a segmentation using $\acv_\mu$ or $\Gcal^n_\mu$ functionals is a good segmentation if it fits properly with the grid of the cells as much as possible, which is the same as resembling the shape of an exact minimizer. In addition, although this does not imply that this fitting behaviour indicates an accurate approximation,  it does make explicit the performance of poor segmentations. To illustrate this application, we consider figure \ref{fig:example-cman}. Here, we segment \ref{fig:cman} via $\acv$ model using two different methods: the anisotropic version of Getreur's implementation (see isotropic one in \cite{Getreur}) and our algorithm \ref{alg:ACVsketch}. Comparing them, we note the
	difference between segmentation \ref{fig:cman-getreur} and  \ref{fig:cman-mp2}. We \sm{observe} that  \ref{fig:cman-getreur} is defective due to the lack of adjustment with the grid of image \ref{fig:cman}, a mismatch that may be rectified by defining other initial parameters. Moreover, this inaccuracy is identified without comparing the segmentation with the exact minimizer \ref{fig:cman-orig}.\\

	In consequence, we consider this application as a \textit{setting up} test which detects incorrect parameter initializations without comparing the approach with the exact minimizer directly.  {In the same spirit, although computing the exact minimizer of $\acv_\mu$ by brute force (i.e. minimizing the functional within the sets in $\mathcal F(G_f)$) is a lengthy computation, it can be used as a test for comparing the performance of any segmentation algorithm in some simple images, in the sense that the grid associated to the image is sufficiently small.}\\

	\subsubsection*{Acknowledgments}
	{\small The authors have been partially supported by
	Conselleria d'Innovació, Universitats, Ciència i Societat Digital, project AICO/2021/223. The first author has also been supported by the Spanish MCIU and FEDER project 094775-B-100. The second author has also been supported by the Universitat de València, grant ``Ajudes per a la col·laboració en la investigació"\ and by the Agència Valenciana de la Innovació, programme ``Promoció del talent'' (Innodocto, Ref. INNTA3/2022/16).}

\bibliographystyle{elsarticle-num} 


\end{document}